\documentclass[leqno,12pt]{article}
\usepackage{latexsym}
\usepackage{graphicx}
\usepackage{amsmath}
\usepackage{amssymb}
\usepackage{mathtools}
\usepackage{cite}
\usepackage{color}
\usepackage{subfig}
\usepackage[margin=1.3in, top=1.3in, bottom=1.3in]{geometry}
\usepackage{algorithm}
\usepackage{algorithmic}
\usepackage{hyperref}
\usepackage{appendix}

\newcommand*{\Nwarrow}{\rotatebox[origin=c]{135}{\(\Longleftrightarrow\)}}
\newcommand*{\Swarrow}{\rotatebox[origin=c]{45}{\(\Longleftarrow\)}}
\newcommand{\nc}{\newcommand}
\nc{\nt}{\newtheorem}
\nt{thm}{Theorem}[section]
\nt{cor}[thm]{Corollary}
\nt{prop}[thm]{Proposition}
\nt{lem}[thm]{Lemma}
\nt{defn}[thm]{Definition}
\nt{rem}[thm]{Remark}
\nt{exa}[thm]{Example}
\nt{ass}[thm]{Assumption}
\nt{alg}[thm]{Algorithm}
\nt{conj}[thm]{Conjecture}
\nt{claim}[thm]{Claim}
\nt{oracle}[thm]{Oracle}
\nc{\ip}[2]{\mbox{$\langle #1,#2 \rangle$}}
\nc{\pf}{\noindent{\bf Proof\ \ }}
\nc{\finpf}{\hfill{$\Box$}\linespace}
\nc{\linespace}{\vspace
{\baselineskip} \noindent}
\nc{\R}{{\mathbf R}}
\nc{\X}{{\mathbf X}}
\nc{\Y}{{\mathbf Y}}
\nc{\E}{{\mathbf E}}
\nc{\oR}{\overline{\R}}
\nc{\M}{\mathcal M}
\nc{\e}{\epsilon}

\nc{\Rn}{{\mathbf R}^n}
\nc{\inT}{\mbox{\rm int}\,}
\nc{\cl}{\mbox{\rm cl}\,}

\def\tto{\;{\lower 1pt \hbox{$\rightarrow$}}\kern -12pt
           \hbox{\raise 2.8pt \hbox{$\rightarrow$}}\;}
\newenvironment{myequation}{\setcounter{equation}{\value{thm}}
   \begin{equation}}{\addtocounter{thm}{1}\end{equation}}

\nc{\bmye}{\begin{myequation}}
\nc{\emye}{\end{myequation}}

\begin{document}
\title{
Recognizing weighted means in geodesic spaces
}
\author{Ariel Goodwin
\thanks{CAM, Cornell University, Ithaca, NY.
\texttt{arielgoodwin.github.io} 
}
\and
Adrian S. Lewis
\thanks{ORIE, Cornell University, Ithaca, NY.
\texttt{people.orie.cornell.edu/aslewis} 
\hspace{2cm} \mbox{~}
Research supported in part by National Science Foundation Grant DMS-2006990.}
\and
Genaro L\'opez-Acedo
\thanks{Department of Mathematical Analysis -- IMUS, University of Seville, 41012 Seville, Spain
\texttt{glopez@us.es} 
}
\and
Adriana Nicolae\thanks{Department of Mathematics, Babe\c{s}-Bolyai University, 400084, Cluj-Napoca, Romania \hfill \mbox{}
\texttt{anicolae@math.ubbcluj.ro} 
}
}
\date{\today}
\maketitle

\begin{abstract}
Geodesic metric spaces support a variety of averaging constructions for given finite sets.  Computing such averages has generated extensive interest in diverse disciplines.  Here we consider the inverse problem of recognizing computationally whether or not a given point is such an average, exactly or approximately.  In nonpositively curved spaces, several averaging notions, including the usual weighted barycenter, produce the same ``mean set''.  In such spaces, at points where the tangent cone is a Euclidean space, the recognition problem reduces to Euclidean projection onto a polytope.  Hadamard manifolds comprise one example.  Another consists of CAT(0) cubical complexes, at relative-interior points:  the recognition problem is harder for general points, but we present an efficient semidefinite-programming-based algorithm. 
\end{abstract}
\medskip

\noindent{\bf Key words:} geodesic space, CAT(0), mean, barycenter, cubical complex, semidefinite program
\medskip

\noindent{\bf AMS Subject Classification:} 90C48, 57Z25, 65K10, 49M29

\section{Introduction}
Given a nonempty finite set $A$ in a metric space $(\X,d)$, we consider the following {\em recognition question\/}:  decide whether or not a given point in $\X$ is some kind of weighted average of the points in $A$.  Such questions are familiar in Euclidean space, but also arise naturally in nonlinear metric spaces:  for example, we might ask whether some phylogenetic tree, viewed as a point in the BHV metric space introduced in \cite{billera}, is an average of several other trees.

\begin{figure}
\centering
\parbox{5cm}
{
\begin{center}
\includegraphics[width=40mm]{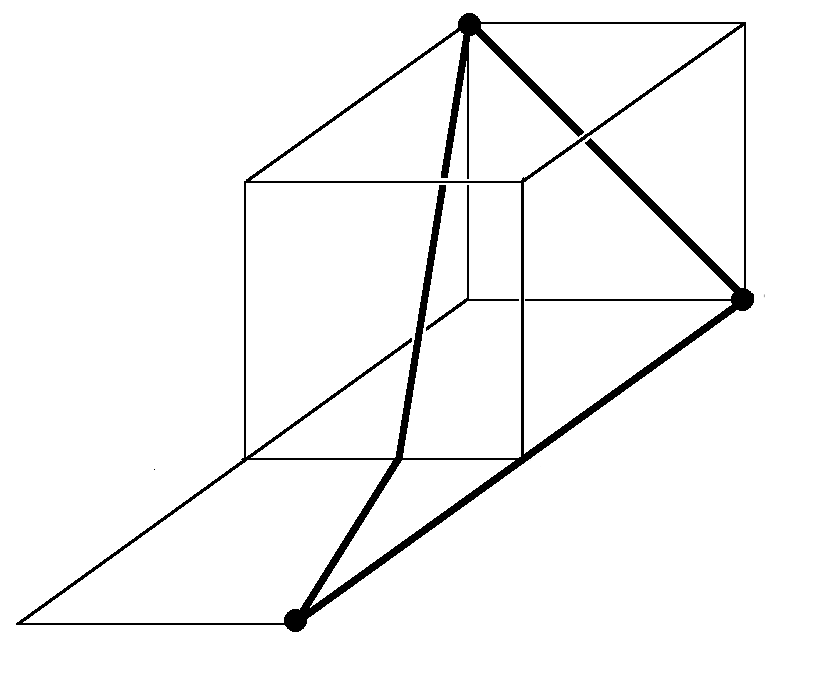}
\end{center}
}
\hspace{20mm}
\parbox{5cm}
{
\begin{center}
\includegraphics[width=40mm]{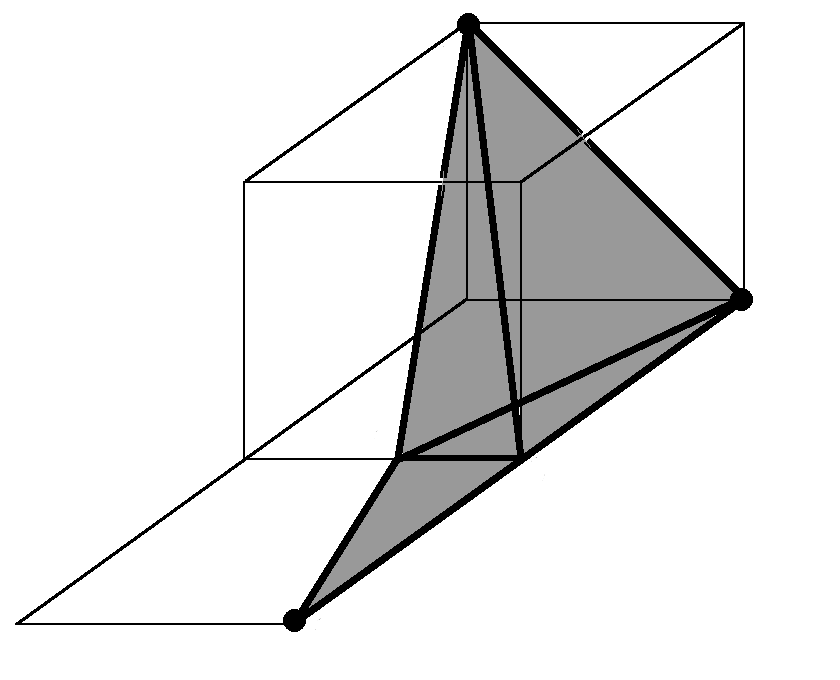}
\end{center}
}
\caption{A triangle and convex hull from three points in a geodesic metric space}
\label{fig:triangle-hull}
\end{figure}

In a vector space, the most immediate notion of a ``weighted average'' is a convex combination.  Recognizing membership in the set of convex combinations from the set $A$ ---  the convex hull of $A$ --- amounts to solving a simple linear program.  The convex hull notion extends naturally, at least to metric spaces that are geodesic.  Figure \ref{fig:triangle-hull} illustrates an example.  The space consists of a three-dimensional cube and a two-dimensional square, each with their Euclidean metrics, glued together along a common edge.  In this space, we consider three points, the geodesics between them, illustrated in the first picture, and their convex hull --- the union of a triangle in the square and a simplex in the cube --- illustrated in the second picture.

In some simple nonlinear spaces, convex hulls are easy to understand:  in any metric tree \cite{hansen}, for example, the convex hull of a finite set just consists of the paths between the various pairs of points.  However, as Figure \ref{fig:triangle-hull} illustrates, in general nonlinear geodesic spaces the convex hull operation behaves less predictably than in Euclidean space.  The convex hull of three points
typically fails to be closed in Riemannian manifolds \cite[Corollary 1.2]{lytchak-generic}, and may fail to be two-dimensional, both in Hadamard manifolds \cite{bessenyi}) and in CAT(0) cubical complexes \cite{lubiw}, as we also see in Figure \ref{fig:triangle-hull}.  Furthermore, in a general Hadamard space, Gromov's question \cite[6.B1(f)]{gromov99} on whether such sets must have compact closure remains open after forty years \cite{bacak-open,basso}.  The recognition question for convex combinations, even in the specific case of a CAT(0) cubical complex,  seems challenging~\cite{lubiw}.

Rather than focusing on convexity, therefore, we consider a simpler version of the recognition question, framed in terms of the most basic metric tools available:  we seek to recognize whether a given point is a barycenter or median.  As well as being easier than recognizing convex combinations, this question may be more interesting in computational practice, reflecting better the notion of a weighted average.

To frame our question more precisely, denote by $\R^A$ the Euclidean space of vectors of components indexed by $A$ with the inner product $\ip{w}{v} = \sum_a w_a v_a$ for vectors $w,v \in \R^A$.  We write $w \ge v$ to mean $w_a \ge v_a$ for all points $a \in A$, and $w > v$ means $w_a > v_a$ for all $a \in A$.  Let 
$\R^A_+$ denote the cone of vectors in $\R^A$ with all nonnnegative components.  Define
the distance vector $d_A \colon \X \to \R^A$ by
\[
\big(d_A(x)\big)_a = d_a(x) = d(a,x) \qquad \mbox{for points}~ x \in \X~ \mbox{and}~ a \in A,
\]
along with its componentwise $p$-power version $d_A^p$, for exponents $p \ge 1$.  We focus on the exponents $p=1,2$ in particular, leading to the notions of {\em weighted medians} and {\em barycenters}, which are minimizers of functions of the form $\ip{w}{d_A}$ and 
$\ip{w}{d^2_A}$ respectively, for nonzero weight vectors $w \in \R^A_+$.  We compare those familiar points with various kinds of Pareto optimal points for the vector objective $d_A$, in particular those points with the property that no other point is strictly closer to every point in $A$.  When the underlying space $\X$ is Euclidean, it is easy to verify that the resulting sets of points all coincide with the convex hull of $A$.
 
\begin{figure}
\centering
\includegraphics[width=50mm]{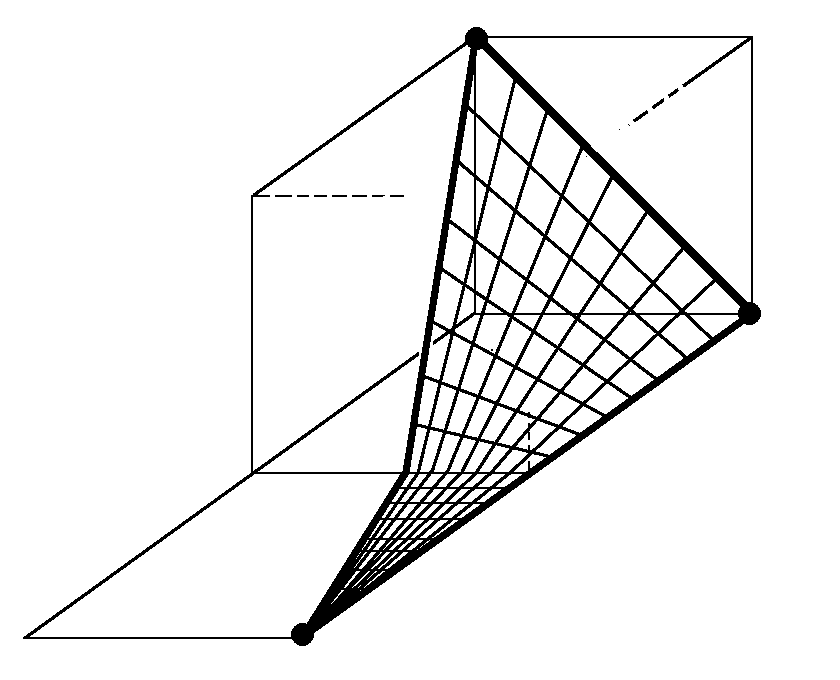}
\caption{The mean set generated by three points in a geodesic metric space}
\label{fig:mean-set}
\end{figure}

In more general geodesic metric spaces, weighted medians and barycenters of a finite set $A$, along with various kinds of Pareto points for the distance vector $d_A$, are easier to understand than the convex hull.  In rather general circumstances, we prove that these sets of points all coincide, comprising a single ``mean set''.  This is the case in particular in Hadamard spaces, where the mean set is a nonempty, compact, path-connected subset of the closed convex hull of $A$.  For the example in Figure \ref{fig:triangle-hull}, we illustrate the mean set in Figure \ref{fig:mean-set}.  

In a Hadamard space $\X$, we prove that a point $\bar x$ belongs to the mean set for a finite set $A$ if and only if it minimizes the ``test'' function 
\[
x ~\mapsto~ \frac{1}{2}\max_{a \in A} \{d_a^2(x) - d_a^2(\bar x) \}.
\]
This function is convex, so its minimizers are points where the metric slope is zero.  We describe the slope in terms of the angle-based geometry of the directions of the geodesics from $\bar x$ towards points in $A$.  In particular, when the tangent cone to $\X$ at $\bar x$ is isometric to a Euclidean space, calculating the slope reduces to a simple projection problem in this space.  Hadamard manifolds have this property, for example, as do CAT(0) cubical complexes at points $\bar x$ in the relative interior of cells that are maximal, such as the cube and the square in Figures \ref{fig:triangle-hull} and \ref{fig:mean-set}.  General points in CAT(0) cubical complexes are more challenging:  we describe a semidefinite-programming-based algorithm to recognize and certify whether or not such points belong to the mean set. 

\section{Means}
Weighted medians and barycenters of a finite set are {\em scalarized} notions of the distance to points in the set, depending on a given weight vector.  Before those ideas, we first consider three natural unscalarized notions that capture the idea of balancing all the points in the set.  

\begin{defn}
{\rm
Given a nonempty finite set $A$ in a metric space $(\X,d)$, we define the following three successively weaker properties of a point $\bar x \in \X$.  
\begin{description}
\item[{\em Strongly $A$-minimal\/}:]
No point $x \ne \bar x$ satisfies $d_A(x) \le d_A(\bar x)$.
\item[{\em $A$-Pareto\/}:]
All points $x \in \X$ such that $d_A(x) \le d_A(\bar x)$ satisfy $d_A(x) = d_A(\bar x)$.
\item[{\em $A$-minimal\/}:]
No point $x \in \X$ satisfies $d_A(x) < d_A(\bar x)$.
\end{description}
We denote the set of $A$-minimal points by $\mbox{minimal}(A)$.  We say that a point $x \in \X$ {\em certifies} the property $\bar x \not\in \mbox{minimal}(A)$ if 
\bmye \label{certificate2}
d_a(x) < d_a(\bar x) \qquad \mbox{for all}~ a \in A.
\emye
The {\em test function} $f_{A,\bar x} \colon \X \to \R$ is defined by
\bmye \label{objective}
f_{A,\bar x}(x) ~=~ \frac{1}{2} \max_{a \in A} \{ d_a^2(x) - d_a^2(\bar x) \} \qquad \mbox{for}~ x \in \X.
\emye
}
\end{defn}

The following result is immediate.

\begin{prop}[Minimizers of the test function] \label{recognition}
Given a nonempty finite set $A$ in a metric space $\X$, a point $\bar x \in \X$ is $A$-minimal if and only if $\bar x$ minimizes the test function $f_{A,\bar x}$, and $\bar x$ is strongly $A$-minimal if and only if $\bar x$ is the unique minimizer of $f_{A,\bar x}$.
\end{prop}

Among scalarized ideas, the following definition covers medians and barycenters.

\begin{defn}
{\rm
Given a nonempty finite set $A$ in a metric space $(\X,d)$, denote the vector in $\R^A$ with all components equal to one by $e$, and define the {\em simplex}
\[
\Delta^A ~=~ \{ w \in \R^A_+ : \ip{e}{w} = 1\}.
\]
For any exponent $p \ge 1$, a point $\bar x \in \X$ is a {\em $p$-mean} for $A$ if there exists a weight vector $w \in\Delta^A$ such that $\bar x$ minimizes $\ip{w}{d_A^p}$.  A $1$-mean is called a {\em median}, and a $2$-mean is called a {\em barycenter}.  The {\em variance function} is $v_A \colon \Delta^A \to \R_+$ defined by 
\[
v_A(w) ~=~ \inf_{x \in \X} \ip{w}{d_A^2(x)} \qquad (w \in \Delta^A).
\]
The {\em mean set} $\mbox{mean}(A)$ is the set of barycenters.  A vector $w$ {\em certifies} the property 
$\bar x \in \mbox{\rm mean}(A)$ when
\bmye \label{certificate1}
w \in \Delta^A \qquad \mbox{and} \qquad \bar x~ \mbox{minimizes}~ \ip{w}{d_A^2}.
\emye
}
\end{defn}

Notice that the variance function is concave.  Other terminology for a barycenter includes {\em center of mass}, {\em center of gravity}, and {\em Fr\'echet} or {\em Karcher mean}.  Existence results for barycenters in Hadamard manifolds date back almost a century, to Cartan~\cite{cartan}.  The idea in more general Hadamard spaces originated with \cite{jost94}.

The following result is immediate.
\begin{prop} \label{easy}
Consider any nonempty finite set $A$ in a metric space, and any exponent $p \ge 1$.  Then any $p$-mean is $A$-minimal, so in particular we have
\[
\mbox{\rm mean}(A) ~\subset~ \mbox{\rm minimal}(A).
\]
\end{prop}

\begin{defn} \label{meb}
{\rm
A {\em circumcenter} of a finite set $A$ in a metric space $(\X,d)$ is a solution of the {\em minimal enclosing ball problem}
\cite{arnaudon-nielsen}
\[
\inf_{x \in \X} \max \{ d_a(x) : a \in A \}.
\]
}
\end{defn}
In passing, we note the following easy observation.

\begin{prop}
Any circumcenter of a finite set $A$ is $A$-minimal.
\end{prop}

In the argument that follows, we interpret the expression $0^0$ as $1$.

\begin{lem} \label{lem}
Consider exponents $q \ge p \ge 1$.  Then any numbers $\alpha,\beta \ge 0$ satisfy the inequality
\[
\frac{1}{q}(\alpha^q - \beta^q) ~\ge~ \beta^{q-p} \cdot \frac{1}{p}(\alpha^p - \beta^p).
\]
\end{lem}

\pf
The case $\beta = 0$ is immediate, so we assume $\beta > 0$.
Since the function $t \mapsto t^{\frac{q}{p}}$ (for $t > 0$) is convex, the subgradient inequality implies
\[
\alpha^q - \beta^q ~=~ (\alpha^p)^{\frac{q}{p}} - (\beta^p)^{\frac{q}{p}}
~\ge~ \frac{q}{p} (\beta^p)^{\frac{q}{p} - 1} (\alpha^p - \beta^p),
\]
from which the result follows.
\finpf

\begin{prop} \label{same}
Consider exponents $q \ge p \ge 1$.  Then any $p$-mean of a finite set $A$ is also a \mbox{$q$-mean}.
\end{prop}

\pf
Consider a $p$-mean $y$.  By definition, there exists a weight vector 
$w \in \Delta^A$ such that $y$ minimizes the function $\ip{w}{d_A^p}$.  By Lemma \ref{lem}, all points $x \in \X$ satisfy 
\[
\frac{1}{q}\big(d_a^q(x) - d_a^q(y)\big) ~\ge~ \big(d_a(y)\big)^{q-p} \cdot 
\frac{1}{p}\big(d_a^p(x) - d_a^p(y)\big).
\]
We deduce
\[
\sum_{a \in A} w_a \big( d_a(y)\big)^{p-q} \big(d_a^q(x) - d_a^q(y)\big) ~\ge~ 0,
\]
so the result follows.
\finpf

In Hilbert space, the relationships between these various properties is summarized in the following  result, whose proof is a combination of standard arguments.

\begin{prop}[Means in a Hilbert space] \label{euclidean}
Consider a Hilbert space $\E$, a nonempty finite set $A \subset \E$, and a point $\bar x \in \E$.  Then the following properties are equivalent:
\begin{enumerate}
\item[{\rm (i)}] $\bar x$ is a median for $A$
\item[{\rm (ii)}] $\bar x$ is a $p$-mean for $A$, for all exponents $p \ge 1$
\item[{\rm (iii)}] $\bar x$ is a barycenter for $A$
\item[{\rm (iv)}] $\bar x$ is $A$-minimal
\item[{\rm (v)}] $\bar x$ is strongly $A$-minimal
\item[{\rm (vi)}] $\bar x$ is $A$-Pareto
\item[{\rm (vii)}] the set $A_{\bar x} = \{z \in \E : \ip{z}{a - \bar x} > 0~ \forall a \in A \}$ is empty
\item[{\rm (viii)}]  $\bar x$ lies in the convex hull of $A$.
\end{enumerate}
In particular, we have
\[
\mbox{\rm mean}(A) ~=~ \mbox{\rm minimal}(A) ~=~ \mbox{\rm conv}(A).
\]
The test function $f_{A,\bar x}$ has a unique minimizer, $\hat x$, which coincides with the nearest point to $\bar x$ in the convex hull of $A$.  The point $\bar x$ is not $A$-minimal if and only if
$\hat x \ne \bar x$, in which case $\hat x - \bar x \in A_{\bar x}$ and the property $\bar x \not\in \mbox{\rm minimal}(A)$ is certified by every point on the geodesic segment $(\bar x,\hat x]$. 
\end{prop}

\pf
The implications (i) $\Rightarrow$ (ii) $\Rightarrow$ (iii) $\Rightarrow$ (iv) follow from Propositions \ref{same} and \ref{easy}.  If property (iv) holds, then the point $\bar x$ minimizes the test function $f_{A,\bar x}$, which is strictly convex, so in fact $\bar x$ is the unique minimizer:  we have therefore proved \mbox{(iv) $\Rightarrow$ (v)}, and (v) $\Rightarrow$ (vi) is immediate.  Since we know (vi) $\Rightarrow$ (iv), we deduce that properties (iv), (v), and (vi) are all equivalent.

Supposed next that property (vii) fails.  Then some point $y \in \E$ satisfies 
$\ip{y - \bar x}{a - \bar x} > 0$ for all points $a \in A$.  
For all small $t>0$ we have
\[
|(ty + (1-t)\bar x) - a|^2  - |\bar x - a|^2
~=~ 2t\ip{\bar x - a}{y-\bar x} + t^2|y-\bar x|^2 ~<~ 0,
\]
so property (iv) fails, and hence so does property (vi).  We have thus proved \mbox{(vi) $\Rightarrow$ (vii)}.  The separating hyperplane theorem proves (vii) $\Rightarrow$ (viii).  Lastly, suppose property (viii) holds, so there exists a weight vector $w$ in the simplex $\Delta^A$ such that $\sum_a w_a a = \bar x$.  We seek to prove that $\bar x$ is a median.  We can suppose $\bar x \not\in A$, since otherwise the result is immediate.  With this assumption, we have
\[
0 ~=~ \sum_{a \in A} w_a (\bar x - a) ~=~ 
\sum_{a \in A} w_a|\bar x - a| \Big( \frac{1}{|\bar x - a|} (\bar x - a) \Big) 
~=~ \sum_{a \in A} w_a|\bar x - a| \nabla d_a(\bar x) .
\]
Hence $\bar x$ minimizes the convex function $\sum_a w_a|\bar x - a| d_a$, and hence is a median.  This completes the proof of the equivalence of the properties (i) -- (viii).

The test function $f_{A,\bar x}$ is strictly convex and weakly lower semicontinuous, with bounded level sets, so has a unique minimizer, $\hat x$.  Consider the value
\[
2f_{A,\bar x}(\hat x) ~=~ \max_{a \in A} \{|\hat x-a|^2 - |\bar x - a|^2\}
~=~ |\hat x|^2 - |\bar x|^2 + 2\max_{a \in A} \ip{\bar x - \hat x}{a} .
\]
Denote by $\bar A$ the subset of $A$ attaining the $\max$, and define a number 
\[
\beta ~=~ f_{A,\bar x}(\hat x) + \frac{1}{2}(|\bar x|^2 - |\hat x|^2).
\] 
We then have
\[
\ip{\bar x - \hat x}{a}
\quad \left\{
\begin{array}{ll}
=~ \beta & (a \in \bar A) \\
< ~ \beta & (a \not\in \bar A).
\end{array}
\right.
\]
In particular, we deduce $\ip{\bar x - \hat x}{x} \le \beta$ for all $x \in \mbox{conv}\,A$.
For all points $x$ near $\hat x$ we have
\[
f_{A,\bar x}(x) ~=~ \frac{1}{2}\max_{a \in \bar A} \{|x-a|^2 - |\bar x - a|^2\}.
\]
Since $\hat x$ minimizes $f_{A,\bar x}$, by convexity it also minimizes the function on the right.  The point $\hat x$ is therefore a barycenter for the set $\bar A$, by our earlier argument.  We therefore know that there exists a vector $w \in \Delta^A$ satisfying $w_a = 0$ for all $a \not\in \bar A$ and such that $\hat x$ minimizes the function $\sum_a w_a |\cdot - a|^2$.  We deduce 
$\hat x = \sum_a w_a a$, so $\ip{\bar x - \hat x}{\hat x} = \beta$.
To conclude, we have $\hat x \in \mbox{conv}\,A$ and $\ip{\bar x - \hat x}{x - \hat x} \le 0$ for all $x \in \mbox{conv}\,A$,
so $\hat x$ is the nearest point to $\bar x$ in $ \mbox{conv}\,A$, as claimed.

The point $\bar x$ is not  $A$-minimal exactly when it does not minimize the test function $f_{A,\bar x}$, which  is equivalent to $\hat x \ne \bar x$.  In that case we have, for all points $a \in A$,
\[
\ip{\hat x - \bar x}{a - \bar x} ~=~ \ip{\hat x - \bar x}{(a - \hat x) + (\hat x - \bar x})
~=~ \ip{\hat x - \bar x}{a - \hat x} ~+~ |\hat x - \bar x|^2 ~>~ 0,
\]
and $|\hat x - a| < |\bar x - a|$, so for all $t \in (0,1)$, by the strict convexity of $|\cdot|^2$, we have
\begin{eqnarray*}
|(t\hat x + (1-t)\bar x) ~-~ a|^2 
&=&
|t(\hat x - a) + (1-t)(\bar x - a)|^2 \\
&<& 
t|\hat x - a|^2 ~+~ (1-t)|\bar x - a|^2 ~\le~ |\bar x - a|^2.
\end{eqnarray*}
Consequently, the point $t\hat x + (1-t)\bar x$ certifies $\bar x \not\in \mbox{mean}(A)$, and we have already proved the case $t=1$. 
\finpf

\section{Duality}
Convex optimization duality plays a fundamental role in our development.  To suggest why this might be, consider the minimal enclosing ball problem, Definition \ref{meb}, for a nonempty finite set $A$ in a general metric space $(\X,d)$.  We can write this problem in the form,
\[
\inf_{x \in \X,\, y \in \R} \{y : y \ge d_a(x)~ \forall a \in A \}.
\]
or equivalently in the saddlepoint form
\[
\inf_{x \in \X,\,y \in \R} \sup_{w \in \R^A_+} \Big\{y + \sum_{a \in A} w_a \big( y - d_a(x) \big) \Big\}.
\]
We can then study this problem through its dual, obtained by interchanging the order of the $\inf$ and the $\sup$.

More generally, we consider an analogous optimization problem associated with the $p$-means of a finite set $A \subset \X$.  We make the following assumption.

\begin{ass} \label{assumption}
{\rm
For a Euclidean space $\Y$, we consider a vector $b \in \R^A$, a nonzero vector $c \in \Y$, and a linear map \mbox{$G \colon \Y \to \R^A$} satisfying the following conditions.
\begin{itemize}
\item
{\bf Strict primal feasibility:}
There exists a point $\hat x \in \X$ and a vector $\hat y \in \Y$ such that $G\hat y - d_A^p(\hat x) > b$.
\item
{\bf Dual feasibility:}
There exists a vector $\hat w \in \R^A_+$ such that $G^*\hat w = c$.
\end{itemize}
}
\end{ass}
With this assumption, given an exponent $p \ge 1$ and a vector $b \in \R^A$, we consider the primal optimization problem
\[
\mbox{(P)} \qquad
\left\{
\begin{array}{ll}
\inf				& \ip{c}{y}							\\
\mbox{subject to}	& Gy - d_A^p(x) ~\ge~ b 	\\	
					& x \in \X,~ y  \in \Y.	 		 
\end{array}
\right.
\]
We can rewrite the problem in the form
\[
\inf_{x \in \X \:,\: y \in \Y } ~ \sup_{w \in \R^A_+} 
\Big \{ \ip{c}{y} ~+~ \ip{w}{b-Gy+d_A^p(x)} \Big\}.
\]
Minimax problems of this form have the following elementary property.

\begin{prop}[Weak minimax duality] \label{minimax}
Given sets $Z$ and $W$ and a function $\phi \colon Z \times W \to \R$, define functions
$\underline\phi \colon Z \to \bar\R$ by $\underline\phi(z) = \sup_w \phi(z,w)$, and
\mbox{$\overline\phi \colon W \to \bar\R$} by $\overline\phi(w) = \inf_z \phi(z,w)$.  Then the inequality 
\[
\inf \underline\phi ~\ge~ \sup \overline\phi
\]
holds.  Suppose furthermore that this inequality in fact holds with equality, and that $\bar z$ minimizes $\underline\phi$ and $\bar w$ maximizes $\overline\phi$. Then, $\bar z$ also minimizes the function $\phi(\cdot,\bar w)$, and $\bar w$ also maximizes the function $\phi(\bar z,\cdot)$.
\end{prop}

With this result in mind, we arrive at the dual problem
\[
\sup_{w \in \R^A_+} ~
\inf_{x \in \X \:,\: y \in \Y }
\Big \{ \ip{c}{y} ~+~ \ip{w}{b-Gy+d_A^p(x)} \Big\},
\]
or equivalently,
\[
\mbox{(D)} \qquad
\left\{
\begin{array}{lrcl}
\sup				& \ip{b}{w} 	& + 	& \inf \ip{w}{d_A^p} 	\\
\mbox{subject to}	& G^*w	& = 	& c					\\
					& w				& \in	& \R^A_+.
\end{array}
\right.
\]
We summarize elementary primal-dual relationships in the following result.

\begin{prop}[Weak duality] \label{weak}
The optimal value of the problem {\rm (P)} is no less than the optimal value of the problem {\rm (D)}. 
If the two values are equal, and $(\bar x,\bar y)$ and $\bar w$ are optimal solutions for {\rm (P)} and {\rm (D)} respectively, then the point $\bar x$ is a $p$-mean of the set $A$, with associated weight vector $\frac{1}{\ip{e}{\bar w}}\bar w$.  
\end{prop}

\pf
Consider the function $\phi \colon (\X \times \Y) \times \R^A_+ \to \R$ defined by
\[
\phi\big( (x,y),w \big) ~=~ \ip{c}{y} ~+~ \ip{w}{b-Gy+d_A^p(x)} \qquad 
(x \in \X,~ y \in \Y,~ w \in \R^A_+).
\]
Then for $(x,y) \in \X \times \Y$ we have
\[
\underline\phi(x,y) ~=~ \sup_{w \in \R^A_+} \phi\big( (x,y),w \big)
~=~
\left\{
\begin{array}{ll}
\ip{c}{y} 	& \mbox{if}~ Gy - d_A^p(x) \ge b  \\
+\infty		& \mbox{otherwise},
\end{array}
\right.
\]
so the problem (P) amounts to minimizing the function $\underline\phi$.  On the other hand, for
$w \in \R_A^+$ we have
\[
\overline\phi(w) ~=~ \inf_{(x,y) \in \X \times \Y} \phi\big( (x,y),w \big) ~=~
\left\{
\begin{array}{ll}
\ip{b}{w} + \inf_{x \in \X} \ip{w}{d_A^p(x)} & \mbox{if}~ G^*w = c \\
-\infty & \mbox{otherwise},
\end{array}
\right.
\]
so the problem $(D)$ amounts to maximizing the function $\overline\phi$.  The first claim now follows from Proposition \ref{minimax}, since $\inf \underline\phi$ is the value of (P) and $\sup \overline\phi$ is the value of (D).  That result also implies that, if the two values are equal, and $(\bar x,\bar y)$ and $\bar w$ are optimal solutions for (P) and (D) respectively, then $G^* \bar w = c$, and furthermore $(\bar x,\bar y)$ minimizes over $\X \times \Y$ the function 
\[
(x,y) ~\mapsto~ \ip{c}{y} ~+~ \ip{\bar w}{b-Gy+d_A^p(x)} ~=~ \ip{\bar w}{b+d_A^p(x)},
\]
so the final claim follows.
\finpf

Following \cite{takahashi}, the metric space $(\X,d)$ is {\em convex} if, for all points $u,v \in \X$ and numbers $\lambda \in [0,1]$ there exists a point $x \in \X$ satisfying the property
\[
d_a(x) ~\le~ \lambda d_a(u) + (1-\lambda)d_a(v) \qquad \mbox{for all points}~ a \in \X.
\]
Many metric spaces in practice are convex.  In particular, any normed space $\X$ is convex, since we can choose $x = \lambda u + (1-\lambda)v$.  Furthermore, CAT(0) spaces, and, more generally, Busemann spaces \cite{papadopoulos} are convex in this sense.  

\begin{thm}[Strong duality] \label{strong}
Suppose that the metric space $(\X,d)$ is convex.  Suppose that Assumption \ref{assumption} holds.
Then the optimal values of the problems {\em (P)} and {\em (D)} are equal and finite, and {\em (D)} has an optimal solution.  Furthermore, for any optimal solution $(\bar x,\bar y)$ for {\em (P)} and $\bar w$ for (D), the point $\bar x$ is a $p$-mean of the set $A$, with associated weight vector 
$\frac{1}{\ip{e}{\bar w}}\bar w$.
\end{thm}

\pf
Define a function $\psi \colon \R^A \to \bar\R$ by
\[
\psi(u) ~=~ \inf_{y \in \Y,~ x \in \X} \{ \ip{c}{y} : Gy - d_A^p(x) \ge u  \}.
\]
For any vector $u \in \R^A$, if $Gy - d_A^p(x) \ge u$, then
\[
0 ~\le~ \ip{\hat w}{Gy - d_A^p(x) - u} ~\le~ \ip{c}{y} - \ip{\hat w}{u},
\]
so $\psi(u) \ge \ip{\hat w}{u} > -\infty$.  On the other hand, for any vector $u \in \R^A$ near $b$, we have $G(\hat y) - d_A^p(\hat x) > u$, and hence $\psi(u) < + \infty$. 

We claim that $\psi$ is convex.  To see this, consider any vectors $u_i \in \R^A$, numbers 
$r_i>\psi(u_i)$, and weights $\lambda_i > 0$, for $i=1,2$, satisfying 
$\lambda_1 + \lambda_2 = 1$.  Then there exist points $x_i \in \X$ and vectors $y_i \in \Y$ satisfying
\[
Gy_i - d_A^p(x_i) \ge u_i  \qquad \mbox{and} \qquad  r_i > \ip{c}{y_i} \qquad \mbox{for}~ i=1,2.
\]
Since $\X$ is convex, there exists a point $x \in X$ satisfying
\[
d_a(x) ~\le~ \lambda_1 d_a(x_1) + \lambda_2 d_a(x_2) \qquad \mbox{for all points}~ a \in A.
\]
Set $y = \lambda_1 y_1 + \lambda_2 y_2$.  Then
\[
Gy - d_A^p(x) - b (\lambda_1 u_1 + \lambda_2 u_2) ~\ge~ \lambda_1 (Gy_1 - d_A^p(x_1) - u_1) ~+~ \lambda_2 (Gy_2 - d_A^p(x_2) - u_2) ~\ge~ 0
\]
and 
\[
\ip{c}{y} ~=~ \lambda_1 \ip{c}{y_1} + \lambda_2 \ip{c}{y_2} ~<~ \lambda_1 r_1 + \lambda_2 r_2,
\]
so $\psi( \lambda_1 u_1 + \lambda_2 u_2) < \lambda_1 r_1 + \lambda_2 r_2$.
Taking the infimum over $r_1$ and $r_2$ now shows the convexity inequality:
$\psi( \lambda_1 u_1 + \lambda_2 u_2) \le \lambda_1 \psi(u_1) + \lambda_2 \psi(u_2)$.

Since the proper convex function $\psi$ is finite near the vector $b$, it has a subgradient there, which we denote $\bar w$.  By definition we have $\ip{\bar w}{u-b} \le \psi(u) - \psi(b)$ for all $u \in \R^A$.
Equivalently, for all points $x \in \X$ and vectors $y \in \Y$ and $u \in \R^A$ satisfying the inequality $Gy - d_A^p(x) \ge u$, we have $\ip{\bar w}{u-b} + \psi(b) \le \ip{c}{y}$.
Making the change of variable $v  = Gy - d_A^p(x) - u$, we can rewrite this property as
\[
\ip{\bar w}{v} ~\ge~ \ip{\bar w}{Gy - d_A^p(x) - b} - \ip{c}{y} ~+~ \psi(b)
\]
for all $x$, $y$, and $v \in \R^A_+$.  We deduce $\bar w \ge 0$, $G^*\bar w = c$, and 
\[
\ip{b}{\bar w} ~+~ \ip{\bar w}{d_A^p(x)} ~\ge~ \psi(b) \qquad \mbox{for all}~ x \in \X.
\]
Thus the inequality $\ip{b}{\bar w} + \inf\ip{\bar w}{d^A} \ge \psi(b)$ holds, so Proposition \ref{weak} implies both the optimality of the vector $\bar w$ for the problem (D) and the remaining claim.
\finpf

\section{Existence of primal solutions}
We aim to prove, under reasonable conditions, that various notions of mean set coincide.  In our argument, geodesic metric spaces play a central role, so before proceeding we review some standard definitions.  We follow the terminology of \cite{bridson}.

In a metric space $\X$, a {\em geodesic segment} is an image of a compact real interval under an isometry into $\X$. We say that $\X$ is a {\em geodesic space} if every two points are joined by a geodesic segment.  A function $f \colon \X \to \R$ is {\em convex} when its compositions with all such isometries are convex, and a set $C \subset \X$ is {\em convex} if $C$ contains all geodesic segments between points in $C$.  A {\em geodesic triangle} $\Delta$ in $\X$ consists of three points along with the geodesic segments joining the three corresponding pairs.  A {\em comparison triangle} for $\Delta$ is a Euclidean triangle with the same three side lengths as $\Delta$.  A geodesic metric space is called {\em CAT(0)} if, in every geodesic triangle, distances between points on its sides are no larger than corresponding distances in a comparison triangle. A {\em Hadamard space} is a complete CAT(0) space:  examples include Hilbert spaces, and all complete, simply-connected Riemannian manifolds with everywhere nonpositive sectional curvature.  The {\em convex hull} of a set $S$ in a Hadamard space is the intersection of all the convex sets containing $S$, and the {\em closed convex hull} is the closure of the convex hull.

Our proof of the equivalence of various notions of mean set will not rely on the existence of solutions of the primal problem (P).   Nonetheless, that existence question is interesting in its own right, and to study it, we focus on geodesic spaces $(\X,d)$ with two further properties.  The first property is from \cite{espinola}.

\begin{defn}
{\rm
A geodesic space $(\X,d)$ is {\em reflexive} if every descending sequence of nonempty, bounded, closed, convex sets has nonempty intersection. 
}
\end{defn}
Examples include Hadamard spaces and reflexive Banach spaces.  For the second property, see \cite{descombes}.

\begin{defn}
{\rm
Given any geodesic space $(\X,d)$, a {\em conical bicombing} is a map 
$\sigma \colon \X \times \X \times [0,1] \to \X$ such that for each pair of points $(x,y) \in \X \times \X$, the function $\sigma(x,y,\cdot)$ is a constant-speed geodesic from $x$ to $y$ (meaning that there exists an isometry $\gamma \colon [0,d(x,y)] \to \X$ satisfying $\gamma(0) = x$, $\gamma\big(d(x,y)\big) = y$, and  $\sigma(x,y,t) = \gamma\big(td(x,y)\big)$ for all $t \in [0,1]$), and the following property holds:
\[
d\big( \sigma(x,y,t) , \sigma(x',y',t) \big) ~\le~ (1-t)d(x,x') + t d(y,y') \qquad \mbox{for all}~ t \in [0,1].
\]
}
\end{defn}

A conical bicombing selects, for each pair of points in the space $\X$, a specific geodesic between them, and guarantees, in particular, the convexity of all the distance functions $d_a$ (for points $a \in \X$) along that geodesic (as can be seen by setting $x'=y'=a$ in the definition).  Thus, in particular, any geodesic space that has a conical bicombing is a convex metric space.  Any normed space has a conical bicombing, namely the map $(x,y,t) \mapsto (1-t)x+ty$.  Busemann spaces (and hence in particular CAT(0) spaces) have conical bicombings:  indeed, a {\em uniquely} geodesic space $(\X,d)$ is Busemann if and only if it has a conical bicombing \cite{papadopoulos}.

The following result concerns reflexive geodesic spaces with a conical bicombing:  examples include all finite-dimensional normed spaces, all uniformly convex Banach spaces, and all Hadamard spaces.

\begin{thm} \label{primal-existence}
Consider a reflexive geodesic space $(X,d)$ with a conical bicombing.  If Assumption~\ref{assumption} holds, then the problem (P) has an optimal solution.
\end{thm}

\pf
By Theorem \ref{strong}, the optimal value $\nu$ of the problem (P) is finite.  For each $k=1,2,3,\ldots$,  define a closed halfspace
\[
H_k ~=~ \Big\{ y \in \Y : \ip{c}{y} \le \nu + \frac{1}{k} \Big\}
\]
and let $F_k$ denote the set of those points $x \in \X$ for which there exists a vector $y \in H_k$ satisfying $d_A^p(x) \le Gy-b$.  By Assumption \ref{assumption}, such $x$ satisfy
\[
\ip{\hat w}{d_A^p(x)} ~\le~ \ip{\hat w}{Gy-b} ~=~ \ip{G^*\hat w}{y} - \ip{\hat w}{b} ~=~ 
\ip{c}{y} - \ip{\hat w}{b} ~\le~ \nu + \frac{1}{k} - \ip{\hat w}{b},
\]
so the sets $F_k$ are all bounded.  Clearly the sets $F_k$ are also nonempty and decreasing. 

The set $G(H_k) - b - \R^A_+$, being a polyhedron, is closed.  Consequently, so is its inverse image under the continuous map $d_A^p$, namely $F_k$.  

For $i=1,2$, given any points $x_i \in F_k$ and weights $\lambda_i \ge 0$ summing to one, there exist vectors $y_i \in \Y$ satisfying $Gy_i - b \ge d_A^p(x_i)$.  The space $\X$ has a conical bicombing, so there exists a constant-speed geodesic $\gamma \colon [0,1] \to \X$ such that $\gamma(0) = x_1$ and 
$\gamma(1) = x_2$, and such that the functions $d_a \circ \gamma$ are convex for all points $a \in A$, and hence so are the functions $d_a^p \circ \gamma$.  Writing $\lambda_1 x_1 + \lambda_2 x_2$ for the point
$\gamma(\lambda_2)$, we deduce
\begin{eqnarray*}
d_A^p(\lambda_1 x_1 + \lambda_2 x_2) & \le & \lambda_1 d_A^p(x_1) + \lambda_2 d_A^p(x_2) 
~\le~ \lambda_1 (Gy_1 - b) + \lambda_2 (Gy_2 - b) \\
& = & G(\lambda_1 y_1 + \lambda_2 y_2) - b,
\end{eqnarray*}
and so $\lambda_1 x_1 + \lambda_2 x_2 \in F_k$.  We have proved that each set $F_k$ is convex.

Since $(\X,d)$ is reflexive, there exists a point $\bar x \in \cap_k F_k$.  Now consider the linear program
\[
\mbox{(P*)} \qquad
\left\{
\begin{array}{ll}
\inf				& \ip{c}{y}							\\
\mbox{subject to}	& Gy - d_A^p(\bar x) ~\ge~ b 	\\	
					& y \in \Y.	 		 
\end{array}
\right.
\]
Clearly the optimal value of \mbox{(P*)} is no smaller than $\nu$.  On the other hand, for each $k=1,2,3,\ldots$, since $\bar x \in F_k$, there exists a vector $y_k \in H_k$ that is feasible for (P*), so the optimal value of (P*) is no larger than
$\nu + \frac{1}{k}$, and hence equals $\nu$.  Let $\bar y$ be any optimal solution of (P*).  Then 
$(\bar x,\bar y)$ is optimal for (P), as required.
\finpf

\section{Examples}
Given a finite set $A$ in a metric space $(\X,d)$, the optimization model (P) covers a variety of interesting problems.  Our main focus is the following special case.

\begin{exa}[intersecting balls] \label{intersecting}
{\rm
Given radii $r_a \ge 0$ corresponding to each point $a \in A$, we ask whether the corresponding closed balls $B_{r_a}(a)$ have empty intersection.  We can express this question via a problem of the form (P).  

Fix any exponent $p \ge 1$.  For the vector $r \in \R^A$ with components $r_a$, interpret the vector \mbox{$r^p \in \R^A$} componentwise.  Then the optimization problem
\[
\mbox{($\mbox{P}_1$)} \qquad
\left\{
\begin{array}{ll}
\inf				& y							\\
\mbox{subject to}	& ye - d_A^p(x) ~\ge~ -r^p 	\\	
					& x \in \X,~ y  \in \R	 		 
\end{array}
\right.
\]
has a feasible solution $(x,y)$ with objective value $y \le 0$ if and only if
\bmye \label{empty}
\bigcap_{a \in A} B_{r_a}(a) \ne \emptyset.
\emye
The problem ($\mbox{P}_1$) arises from problem (P) by choosing the space $\Y = \R$, the vectors $c = 1$ and $b=-r^p$, and defining the map $G$ by $Gy = ye$ for numbers $y \in \R$.  Assumption~\ref{assumption} holds, since we can choose the vector $\hat w = \frac{1}{|A|}e$, select any point $\hat x \in \X$, and then choose any number $\hat y > \max_a \{ d_a^p(\hat x) - r_a^p \}$.

The dual problem is
\[
\mbox{($\mbox{D}_1$)} \qquad \sup_{w \in \Delta^A} \{ \inf\ip{w}{d_A^p} -\ip{w}{r^p} \}.
\]
In the special case $p=2$, we can rewrite using the variance function: 
\[
\sup_{w \in \Delta^A} \{ v_A(w) -\ip{w}{r^2} \}.
\]
Suppose that the metric space $(\X,d)$ is convex.  Then, by Theorem \ref{strong}, the optimal values of the problems ($\mbox{P}_1$) and ($\mbox{D}_1$) are equal and finite, and ($\mbox{D}_1$) has an optimal solution. Furthermore, for any optimal solution $(\bar x,\bar y)$ for ($\mbox{P}_1$) and $\bar w$ for 
($\mbox{D}_1$), the point $\bar x$ is a $p$-mean for the set $A$, with associated weight vector $\bar w$.  

Assuming further that $(\X,d)$ is a reflexive geodesic space with a conical bicombing, or in particular a Hadamard space, then the problem ($\mbox{P}_1$) also must have an optimal solution $(\bar x,\bar y)$.  Any such solution determines the answer to our original question:  $\bar y \le 0$ if and only if equation (\ref{empty}) holds.
}
\end{exa}

Two special cases of Example \ref{intersecting} are especially interesting.  The first is particularly simple.

\begin{exa}[circumcenters]
{\rm
A point $\bar x \in \X$ is a circumcenter for the set $A$ if and only if the pair $\big(\bar x, \max_a d^p_a(\bar x)\big)$ is an optimal solution of the problem ($\mbox{P}_1$) in the case $r=0$.  Consequently, in convex metric spaces, circumcenters must be $p$-means.  In reflexive geodesic spaces with conical bicombings, we also deduce the well-known result that circumcenters must exist \cite[Example 2.2.18]{bacak-convex}.
}
\end{exa}

The second special case of Example \ref{intersecting} is central to our development.

\begin{exa}[$A$-minimal points] \label{pareto}
{\rm
A point $\bar x \in \X$ is $A$-minimal if and only if the pair $(\bar x, 0)$ is an optimal solution of the problem ($\mbox{P}_1$) in the case $r = d_A(\bar x)$.  Consequently, in a convex metric space, $A$-minimal points must be $p$-means.
In the special case $p=2$, the dual problem is
\bmye \label{p=2}
\sup_{w \in \Delta^A} \{ v_A(w) -\ip{w}{d_A^2(\bar x)} \}.
\emye
}
\end{exa}

The model ($\mbox{P}$) also covers some interesting problems not of the form ($\mbox{P}_1$).

\begin{exa}[$p$-means]
{\rm
Given a nonzero weight vector $w^* \in \R^A_+$, a point $\bar x \in \X$ is a corresponding $p$-mean if and only if the pair $\big(\bar x, d^p_A(\bar x)\big)$ is an optimal solution of the problem
\[
\left\{
\begin{array}{ll}
\inf				& \ip{w^*}{y}							\\
\mbox{subject to}	& y - d_A^p(x) ~\ge~ 0 	\\	
					& x \in \X,~ y  \in \R^A.	 		 
\end{array}
\right.
\]
This corresponds to the problem (P) for the space $\Y = \R^A$, the vectors $c = w^*$ and $b=0$, and the map $G$ equal to the identity.  Assumption \ref{assumption} holds, for any point $\hat x \in \X$, any vector 
$\hat y > d_A^p(\hat x)$, and $\hat w = w^*$.  The dual problem is trivial:
\[
\sup_{w=w^*} \inf \ip{w}{d_A^p}.
\]
}
\end{exa}

\begin{exa}[projection onto intersecting balls]
{\rm
Given a point $h \in \X$, a finite set $\bar A \subset \X$, and a vector $r \in \R^{\bar A}_+$ whose components are radii of the closed balls $B_{r_a}(a)$ centered at each point $a \in \bar A$, a point 
$\bar x \in \X$ is a nearest point to the point $h$ in the set
\[
\bigcap_{a \in \bar A} B_{r_a}(a)
\]
if and only if the pair 
$\big(\bar x, d^p_h(\bar x)\big))$ is an optimal solution of the problem
\[
\left\{
\begin{array}{ll}
\inf				& y							\\
\mbox{subject to}	& y - d^p_h(x) ~ \ge	~ 0 	\\	
					& \mbox{\hspace{5pt}}   - d^p_{\bar A}(x) ~\ge~ -r^p \\
					& y  \in \R,~ x \in \X.			 
\end{array}
\right.
\]
This corresponds to the problem (P), for the set $A = \{h\} \cup \bar A$, the space $\Y = \R$, the vectors $c = 1$ and $b=(0,-r^p)$, and the map $G$ defined by $Gy = (y,0) \in \R \times \R^{\bar A}$.  Assumption \ref{assumption} holds, providing that there exists a point $\hat x$ satisfying $d_{\bar A}(\hat x) < r$, by setting any number $\hat y > d_h^p(\hat x)$, and choosing $\hat w = (1,0) \in \R \times \R^{\bar A}$.
}
\end{exa}

\section{Means in Hadamard spaces}
Our principal focus in this work is to compare the various notions of weighted averages that we have discussed in the particular case of Hadamard spaces.  We begin with the following observation.

\begin{thm} \label{coincide}
Consider a Hadamard space $\X$ (or, more generally, a convex metric space).
For any nonempty finite set $A \subset \X$, the set of medians, barycenters, and more generally of $p$-means (for any exponent $p \ge 1$) all coincide with the set of $A$-minimal points.  In particular, we have
\[
\mbox{\rm mean}(A) ~=~ \mbox{\rm minimal}(A).
\]
\end{thm}

\pf
This follows from Proposition \ref{easy}, Theorem \ref{strong}, and Example \ref{pareto}.
\finpf

\begin{defn}
{\rm
For any finite subset $A$ of a metric space $\X$, the {\em barycenter map} 
$\Phi^A \colon \Delta^A \tto \X$ is defined by
\[
\Phi^A(w) ~=~ \mbox{argmin}\ip{w}{d_A^2} \qquad (w \in \Delta^A).
\]
}
\end{defn}
The following theorem collects facts from \cite[Propositions 4.3, 4.4, and 6.1, and Theorem~6.3]{sturm2}.

\begin{thm} \label{sturm}
For any finite subset $A$ of a Hadamard space, the barycenter map $\Phi^A$ is single-valued and continuous on the simplex $\Delta^A$, with range contained in the closed convex hull of $A$.    For any weight vector $w \in \Delta^A$, the corresponding barycenter $\bar x = \Phi^A(w)$ satisfies the inequality
\bmye \label{variance}
\ip{w}{d_A^2(x) - d_A^2(\bar x)} ~\ge~ d^2(x,\bar x) \qquad \mbox{for all points}~ x \in \X.
\emye
\end{thm}

\begin{cor} \label{topology}
Consider any finite subset $A$ of a Hadamard space. The mean set is a nonempty, compact, path-connected subset of the closed convex hull of $A$, and it coincides with the set of strongly $A$-minimal points and the set of $A$-Pareto points.  In particular, we have
\[
\mbox{\rm mean}(A) ~=~ \mbox{\rm minimal}(A) ~\subset~ \mbox{\rm cl}\, \mbox{\rm conv}(A).
\]
\end{cor}

\pf
Consider the following four properties of a point $\bar x$:
\begin{enumerate}
\item[(i)]
$\bar x \in \mbox{\rm mean}(A)$
\item[(ii)]
$\bar x$ is strongly $A$-minimal
\item[(iii)]
$\bar x$ is $A$-Pareto
\item[(iv)]
$\bar x \in \mbox{\rm minimal}(A)$
\end{enumerate}
Property (i) states that there exists a weight vector $w \in \Delta^A$ satisfying $\bar x = \Phi^A(w)$.  If property (ii) fails, then there exists a point $x \ne \bar x$ satisfying $d_A(x) \le d_A(\bar x)$, and hence $\ip{w}{d_A^2(x) - d_A^2(\bar x)} \le 0$, contradicting inequality (\ref{variance}).  We have proved the implication (i) $\Rightarrow$ (ii).  The implications (ii) $\Rightarrow$ (iii) $\Rightarrow$ (iv) are immediate, and the implication (iv) $\Rightarrow$ (i) follows from Theorem \ref{coincide}.  The four properties are therefore equivalent.  The remainder of the result now follows from  Theorem \ref{sturm}, since $\Delta^A$ is nonempty, compact and path-connected.
\finpf

\noindent
We deduce a standard property of Hadamard spaces.

\begin{cor} \label{two}
The mean set of any two points in a Hadamard space coincides with their convex hull, namely the geodesic segment joining them.
\end{cor}

Comparing Corollary \ref{topology} with the special case of Hilbert space described in Proposition \ref{euclidean}, we see two significant distinctions.  First, unlike in Hilbert space, Corollary \ref{topology} only claims a one-sided relationship between mean sets and closed convex hulls.  Indeed, although equality holds for sets of two points (Corollary \ref{two}), three points in a Hadamard space may have a mean set that is much smaller than their convex hull.  In a CAT(0) cubical complex, \cite{lubiw} presents an example of a set $A$ of three points whose convex hull contains a three-dimensional ball, and indeed we discuss a similar example in Section \ref{mean-example}.  On the other hand, for any set $A$ of cardinality $k$, since the simplex $\Delta^A$ is a $(k-1)$-dimensional semi-algebraic set, and the barycenter map $\Phi^A$ in the proof of Theorem \ref{topology} is semi-algebraic, the mean set, which is the range $\Phi^A(\Delta^A)$, can have dimension at most $k-1$.  The example in Section \ref{mean-example} illustrates this behavior in the case $k=3$.  A similar argument applied to an example of \cite{bessenyi} shows the same phenomenon in a Hadamard manifold.

The second distinction is the absence of any analogue in Corollary \ref{topology} of the angle condition, property (vii), in Proposition \ref{euclidean}. We remedy that as follows.  We refer to {\em angles} in the sense of Alexandrov \cite[Definition I.1.12]{bridson}.

\begin{defn}
{\rm
Given a finite set $A$ in a Hadamard space $\X$, a point $\bar x \in \X$ is {\em $A$-balanced} if either $\bar x \in A$ or $\bar x \not\in A$ and for every point $y \ne \bar x$ in $\X$ there is point $a \in A$ such that the angle $\angle y\bar x a$ is at least~$\frac{\pi}{2}$.
}
\end{defn}

\noindent
The inclusion, in Proposition \ref{euclidean}, of property (vii) in the list of equivalences is then captured by the following result.

\begin{prop}[Recognizing means by angles]~ \label{angle}
Given a nonempty finite subset $A$ of a Hadamard space $\X$,  a point $\bar x \in \X$ is $A$-minimal if and only if it is $A$-balanced.
\end{prop}

\pf
For points $y \ne \bar x \ne a$ in $\X$, the first variation formula \cite[Corollary~4.5.7]{burago} states
\bmye \label{first}
\frac{1}{t} \Big( d_a \big( (1-t)\bar x + ty \big) - d_a(\bar x) \Big) ~\downarrow~ 
-d(\bar x,y)\cos\angle y\bar x a \qquad \mbox{as}~ t \downarrow 0.
\emye
By definition, for any point $\bar x \not\in \mbox{minimal}(A)$, there exists a point $y \in \X$ satisfying 
$d_a(y) < d_a(\bar x)$ for all points $a \in A$, so the property (\ref{first}) shows 
$\angle y\bar x a < \frac{\pi}{2}$, completing the proof in one direction.

Conversely, consider any point $\bar x \in \mbox{minimal}(A)$ and any point $y \ne \bar x$ in $\X$.  Consider any sequence $(t_r)$ in the interval $(0,1)$ decreasing to $0$.  By definition, for each $r=1,2,3,\ldots$, there exists a point $a_r \in A$ satisfying 
$d_{a_r} \big( (1-t_r)\bar x + t_ry \big) \ge d_{a_r}(\bar x)$.
Since $A$ is finite, after taking a subsequence we can suppose that there is a point $a \in A$ such that
$d_a \big( (1-t_r)\bar x + t_ry \big) \ge d_a(\bar x)$ for all $r=1,2,3,\ldots$.
Formula (\ref{first}) now implies $\angle y\bar x a \ge \frac{\pi}{2}$, as required.
\finpf

Given a finite set $A$ in a metric space $\X$, we can summarize the relationship between the various properties that we have defined for a point in $\X$ as follows.
\medskip

\begin{eqnarray*}
& \mbox{median of $A$} & \\ \\
\Swarrow & & \Nwarrow \quad \mbox{{\em (if $\X$ is convex)}} \\ \\
\mbox{$A$-balanced} \quad \Longleftrightarrow \quad \mbox{barycenter of $A$} & \Longrightarrow & \mbox{$A$-minimal} \\ \\
\big\Downarrow  \qquad \mbox{{\em (if $\X$ is Hadamard)}} \qquad \big\Updownarrow \qquad \mbox{} & & \mbox{} \qquad \big\Uparrow \\ \\
\mbox{in $\mbox{cl}\,\mbox{conv}(A)$} \quad \Longleftarrow \quad \mbox{strongly $A$-minimal} & \Longrightarrow & \mbox{$A$-Pareto}
\end{eqnarray*}
\medskip

\section{The mean deficit} \label{mean-deficit}
We introduce the following definition.  For any value $\alpha \in \R$, we denote the positive part 
$\max\{\alpha,0\}$ by $\alpha^+$.

\begin{defn} \label{deficit-def}
{\rm
Given a nonempty finite set $A$ in a Hadamard space, the {\em mean deficit} function on the complement $A^c$ is defined by
\bmye \label{deficit}
\mbox{def}_A(x) ~=~ \Big( \sup_{y \ne x} \: \min_{a \in A} \{ d_a(x)\cos\angle axy \} \Big)^+
\qquad (x \not\in A).
\emye
}
\end{defn}

\noindent
The positive part is important to rule out negative values.  Consider, for example, the case
$A = \{1,-1\} \subset \R = \X$ and the point $x=0$.

\begin{prop}[Euclidean mean deficit] \label{euclid-deficit}
The mean deficit for a nonempty finite set $A$ in a Euclidean space is given by
\[
\mbox{\rm def}_A(x) ~=~ \mbox{\rm dist}(x , \mbox{\rm conv}\, A)
\qquad (x \not\in A).
\]
\end{prop}

\pf
If $\mbox{\rm def}_A(x) > 0$, then the von Neumann minimax theorem implies
\begin{eqnarray*}
0 ~<~ \mbox{\rm def}_A(x) 
&=& 
\sup_{y \ne x} \: \min_{a \in A} \{ |a-x|\cos\angle axy \} ~=~
\sup_{|y - x|=1} \: \min_{a \in A} \{ |a-x|\cos\angle axy \} \\
&=&
\sup_{|y - x|=1} \: \min_{a \in A} \ip{a-x}{y-x} ~=~
\max_{|u|=1} \: \min_{a \in A} \ip{a-x}{u} \\
&=&
\max_{|u| \le 1} \: \min_{a \in A} \ip{a-x}{u} ~=~
\max_{|u| \le 1} \: \min_{a \in \mbox{\scriptsize conv}\,A} \ip{a-x}{u} \\
&=& 
\min_{a \in \mbox{\scriptsize conv}\,A} \max_{|u| \le 1} \ip{a-x}{u}
~=~
\min_{a \in \mbox{\scriptsize conv}\,A} |a-x|,
\end{eqnarray*}
as required.  On the other hand, if $\mbox{\rm def}_A(x) = 0$, then $x \in \mbox{conv}\, A$, by Proposition~\ref{euclidean}, completing the proof.
\finpf

\noindent
With this terminology, we can restate Proposition \ref{angle} as follows.

\begin{cor} \label{coincide2}
Given a nonempty finite subset $A$ of a Hadamard space $\X$,  a point $x \not\in A$ is $A$-minimal if and only if $\mbox{\rm def}_A(x) = 0$.
\end{cor}

Given a point $\bar x \in \X$, we again consider the test function:  the convex function $f_{A,\bar x} \colon \X \to \R$ defined by
\[
f_{A,\bar x}(x) ~=~ \frac{1}{2} \max_{a \in A} \{ d^2_a(x) - d^2_a(\bar x) \} \qquad (x \in \X).
\]
The test function and the mean deficit are closely related through the idea of the metric slope.  
For any metric space $(\X,d)$, function $f \colon \X \to \bar\R$, and point $\bar x \in \X$ where the value $f(\bar x)$ is finite, the {\em slope} $|\nabla f|(\bar x)$ is the quantity
\[
\limsup_{x \to \bar x} \frac{f(\bar x) - f(x)}{d(\bar x,x)},
\]
unless $\bar x$ minimizes $f$ locally, in which case the slope is zero.  When $f$ is convex, the slope is zero if and only if $\bar x$ minimizes $f$, so as a consequence of Theorem~\ref{coincide} and Proposition \ref{recognition}, we have the following result.

\begin{cor}
Consider a Hadamard space $\X$ and a nonempty finite set $A \subset \X$.  A point $\bar x \in \X$ lies in the mean set for $A$ if and only if the slope of the test function satisfies 
$|\nabla f_{A,\bar x}|(\bar x) = 0$.
\end{cor}  

As we now show, under reasonable conditions, the slope of the test function is determined by the mean deficit.

\begin{thm}[Test function slope]
Consider a Hadamard space $\X$, a nonempty finite set $A \subset \X$, and a point $\bar x \not\in A$.  The slope of the test function always satisfies the inequality
\bmye \label{slope-bigger}
|\nabla f_{A,\bar x}|(\bar x) ~\ge~ \mbox{\rm def}_A(\bar x).
\emye
If $\X$ is locally compact, with the geodesic extension property, then equality holds, and furthermore the supremum in definition (\ref{deficit}) of the mean deficit $ \mbox{\rm def}_A(\bar x)$ is attained.
\end{thm}

\pf
To simplify notation, we write $f$ for the test function $f_{A,\bar x}$.
Consider any point $y \ne \bar x$.  For $t \in [0,1]$, denote the convex combination $(1-t)\bar x + ty$ by $x_t$.

The slope satisfies
\[
|\nabla f|(\bar x) ~\ge~ 
\limsup_{1 \ge t \downarrow 0} 
\frac{f(\bar x) - f(x_t)}{d(\bar x,x_t)}
~=~
\limsup_{1 \ge t \downarrow 0} 
\frac{-f(x_t)}{td(\bar x,y)}
\]
Choose any sequence $t_r \in (0,1)$ decreasing to $0$, and any corresponding sequence of points 
$a_r \in A$ such that $f(x_{t_r}) = d^2_{a_r}(x_{t_r}) - d^2_{a_r}(\bar x)$ for all $r$.  After taking a subsequence, we can suppose there exists a point $a \in A$ satisfying $a_r = a$ for all $r$.  We deduce, using the first variation formula, the inequalities
\begin{eqnarray*}
|\nabla f|(\bar x) 
&\ge& 
\limsup_{r \to \infty} \frac{-f(x_{t_r})}{t_rd(\bar x,y)}
~=~
\frac{1}{2} \limsup_{r \to \infty} \frac{d^2_{a}(\bar x) - d^2_{a}(x_{t_r})}{t_rd(\bar x,y)} \\
&=& 
\frac{1}{2} \limsup_{r \to \infty} 
\frac{\big(d_{a}(\bar x) - d_{a}(x_{t_r})\big)\big(d_{a}(\bar x) + d_{a}(x_{t_r})\big)}{t_rd(\bar x,y)} 
\\
&=&
d_a(\bar x) \cos\angle a\bar x y ~ \ge ~
\min_{a \in A} \{d_a(\bar x) \cos\angle a\bar x y \}.
\end{eqnarray*}
Since $y \ne \bar x$ was arbitrary, we have proved
\[
|\nabla f|(\bar x) ~\ge~ \sup_{y \ne \bar x} \: \min_{a \in A} \{ d_a(\bar x)\cos\angle a\bar xy \}.
\]
Since $|\nabla f|(\bar x) \ge 0$, we have proved inequality (\ref{slope-bigger}).

If the point $\bar x$ minimizes the test function $f$, then the slope $|\nabla f|(\bar x)$ is zero.  On the other hand, Proposition \ref{angle} implies $\mbox{def}_A(\bar x) = 0$, so inequality (\ref{slope-bigger}) holds with equality.  Assume now, therefore, that $\bar x$ does not minimize the test function $f$.  By convexity, $\bar x$ is not a local minimizer of $f$.   Therefore there exists a sequence of points $x_r \to \bar x$ in $\X$ satisfying $0 < d(\bar x,x_r) < 1$ for each $r=1,2,3,\ldots$, and
\[
\frac{f(\bar x) - f(x_r)}{d(\bar x,x_r)} ~\to~ |\nabla f|(\bar x).
\]

Henceforth, suppose that the Hadamard space $\X$ is locally compact, and has the geodesic extension property.  For each $r=1,2,3,\ldots$, there exists a point $y_r \in \X$ satisfying $d(\bar x,y_r) = 1$ such that setting $t_r = d(\bar x,x_r)$ gives $x_r = (1-t_r)\bar x + t_r y_r$.  Since $t_r \to 0$, we therefore have
\[
\frac{1}{t_r}\min_{a \in A} \{ d^2_a(\bar x) - d^2_a(x_r)  \}
~\to~ |\nabla f|(\bar x).
\]
By local compactness, after taking a subsequence, we can suppose that the points $y_r$ converge to some point $y \in \X$, which must therefore satisfy $d(\bar x, y) = 1$.  

Consider any point $a \in A$.  We must have
\[
|\nabla f|(\bar x) ~ \le ~ \frac{1}{2} \liminf_{r \to \infty} 
\frac{1}{t_r} \{ d^2_a(\bar x) - d^2_a(x_r)  \} 
~ = ~
d_a(\bar x) \liminf_{r \to \infty} 
\frac{1}{t_r} \{ d_a(\bar x) - d_a(x_r)  \}.
\]
Define points $x'_r = (1-t_r)\bar x + t_r y$ for each $r=1,2,3,\ldots$.
Since the distance function $d_a$ is $1$-Lipschitz, and the metric $d$ is jointly convex \cite[Proposition II.3.3]{bridson}, we have
\[
\frac{1}{t_r}| d_a(x_r) ~-~  d_a(x'_r) |
~\le~ \frac{1}{t_r}d ( x_r,x'_r)
~\le~ d(y_r,y) ~\to~ 0,
\]
so using the first variation formula again shows
\[
|\nabla f|(\bar x) ~ \le ~ d_a(\bar x)\liminf_{r \to \infty} 
\frac{1}{t_r} \{ d_a(\bar x) - d_a\big((1-t_r)\bar x + t_r y\big)  \}
~=~  d_a(\bar x)\cos\angle a\bar xy .
\]
Since $a$ was arbitrary, we have proved
\[
|\nabla f|(\bar x) ~\le~ \min_{a \in A} d_a(\bar x)\cos\angle a\bar xy.
\]
We conclude that in fact
\[
|\nabla f|(\bar x) ~=~  \sup_{y \ne \bar x} \min_{a \in A} d_a(\bar x)\cos\angle a\bar xy,
\]
and the supremum is attained.

To prove that the supremum is attained in full generality, we note, using the geodesic extension property, that we can restrict attention in the definition (\ref{deficit}) of the mean deficit to points $y \in \X$ satisfying $d(\bar x,y) =1$.  By local compactness, this set is compact, and the angle 
$\angle a\bar x y$ depends continuously on the point $y$, by \cite[Proposition II.2.2]{bridson}, so attainment follows.
\finpf

Notice that Proposition \ref{euclid-deficit} fails in non-Euclidean settings:  the mean deficit may not equal the distance to either the convex hull or the mean set, and may fail to be continuous, as the following example shows.

\begin{exa} \label{tripod}
{\rm
The {\em tripod} is the Hadamard space consisting of three copies of the half-line $\R_+$ glued together at the common point $0$.  Consider the set $A$ consisting of two points, each on distinct half-lines at unit distance from $0$.  For any value $t>0$, consider the point $x_t$ on the third half-line at a distance $t$ from $0$. Then we have $\mbox{def}_A(x_t) = 1+t$ for all $t > 0$ and $x_t \to 0$, but
$\mbox{def}_A(0) = 0$.
}
\end{exa}

\begin{prop}[Mean deficit level sets]~ \label{bounded2}
Consider a Hadamard space $\X$ and a nonempty finite set $A \subset \X$. The mean deficit ${\rm def}_A$ is lower semicontinuous throughout the complement $A^c$.  For all $\kappa \ge 0$, the set
\bmye \label{deficit-level}
A ~\cup~ \{ x \not\in A : \mbox{\rm def}_A (x) \le \kappa \}
\emye
is closed and bounded.
\end{prop}

\pf
Consider any value $\kappa \ge 0$ and any sequence $x_r \to \bar x \not\in A$ satisfying 
\mbox{$\mbox{\rm def}_A(x_r) \le \kappa$} for all $r=1,2,3,\ldots$. Fix any point $y \ne \bar x$. After taking a subsequence, we can suppose, for each $r$, we have $y \ne x_r$, and hence  
\[
\min_{a \in A} \: d_a(x_r) \cos\angle ax_ry ~\le~ \kappa.
\]
Choose $a_r \in A$ attaining this minimum.  After taking a further subsequence, we can suppose that some point $\bar a \in A$ satisfies $a_r = \bar a$ for all $r$, so
$d_{\bar a}(x_r) \cos\angle \bar a x_ry \le \kappa$ for all $r$.
The function $x \mapsto \angle \bar a xy$ is upper semicontinuous at the point $\bar x$, by \cite[Proposition~II.3.3]{bridson}, so  $x \mapsto \cos\angle \bar a xy$ is lower semicontinuous, and hence we deduce 
\[
\min_{a \in A} \: d_a(\bar x) \cos\angle a \bar x y ~\le~ d_{\bar a}(\bar x) \cos\angle \bar a\bar x y ~\le~ \kappa.
\]
Since $y$ was arbitrary, we deduce $\mbox{\rm def}_A(\bar x) \le \kappa$, proving lower semicontinuity.

We deduce that the set (\ref{deficit-level}) is closed, so it remains only to prove that it is bounded.  Suppose, by way of contradiction, that it contains an unbounded sequence $(x_r)$.  Fix any point $y \in \X$.  After taking a subsequence, we can suppose $y \ne x_r \not\in A$ for all $r=1,2,3,\ldots$.  Choose $a_r \in A$ so that $a=a_r$ attains 
$\min_a d_a(x_r) \cos\angle a x_r y$.  After taking a further subsequence, we can suppose that, for all $r$, some point $a \in A$ satisfies $a_r = a$, and so
$d_a(x_r) \cos\angle ax_ry \le \kappa$.  By the law of cosines, we have
\begin{eqnarray*}
2d(x_r,y)\big(d(x_r,y) - d(a,y)\big) & = & \big(d(x_r,y) - d(a,y)\big)^2 + d^2(x_r,y) - d^2(a,y) \\
& \le & d^2(x_r,a) + d^2(x_r,y) - d^2(a,y) ~\le ~ 2\kappa d(x_r,y).
\end{eqnarray*}
We deduce $d(x_r,y) - d(a,y) \le~ \kappa$,
contradicting the unboundedness of $(x_r)$.
\finpf

The mean deficit satisfies the following proximity property for the mean set.

\begin{prop}[Mean deficit and proximity to the mean set]\hspace{1pt} \label{proximity}
Consider a nonempty finite subset $A$ of a locally compact Hadamard space $\X$.  For any sequence $(x_r)$ in the complement $A^c$, we have 
\[
\mbox{\rm def}_A(x_r) \to 0 \qquad \Rightarrow \qquad \mbox{\rm dist}\big(x_r , \mbox{\rm mean}(A) \big) \to 0.
\]
\end{prop}

\pf
Consider a sequence $(x_r)$ satisfying the first property but not the second.  By Proposition \ref{bounded2}, $(x_r)$ is bounded, and furthermore, there exists a value $\epsilon > 0$ and a subsequence $R$ of the positive integers such that $\mbox{\rm dist}\big(x_r , \mbox{\rm mean}(A) \big) \ge \epsilon$ for all $r \in R$.  By local compactness, there exists a further subsequence $R' \subset R$ such that $x_r$ converges to some limit $\bar x$ as $r \to \infty$ in $R'$.  By lower semicontinuity we deduce $\mbox{\rm def}_A(\bar x) = 0$. Proposition \ref{coincide2} now shows $\bar x \in \mbox{mean}(A)$, so $d(x_r,\bar x) \ge \e$ for all $r \in R'$, which is a contradiction.
\finpf 

\noindent
Example \ref{tripod} shows that the converse implication in Proposition \ref{proximity} can fail.

\subsection*{The mean deficit and the tangent cone}
For a point $x$ in the Hadamard space $\X$, Definition \ref{deficit-def} for the mean deficit associated with a set $A \subset \X$ concerns only the distances $d_a(x)$ for each point $a \in A$ and the angles between the geodesic segments $[x,a]$ and other geodesics from $x$.  It is therefore illuminating to view it in terms of the {\em tangent cone} for $\X$ at $x$, written $T_{x}\X$.  For a summary of the construction of the tangent cone, see \cite[Section II.3]{bridson} and \cite{lewis-lopez-nicolae-basic}.  Elements of the tangent cone are pairs $([\gamma],r)$, where $\gamma$ is a geodesic with origin $x$ and 
$[\gamma]$ denotes the equivalence class of geodesics subtending with $\gamma$ an angle zero at $x$, and $r \ge 0$ is a scalar. On this metric cone is defined the following scalar product:
$\ip{([\gamma],r)}{([\gamma'],r')} = rr'\cos\angle(\gamma,\gamma')$.

We can now rewrite the definition of the mean deficit using the terminology of tangent cones.  The resulting formulation is as follows:
\[
\mbox{def}_A(x) ~=~
\Big( \sup_{v = ([\gamma],1) \in T_x \X} \: \min_{a \in A} \Big< v \: , \: \big(\big[[x,a]\big],d_a(x)\big) \Big> \Big)^+ \qquad (x \not\in A).
\]

This reformulation is just a change of language, emphasizing the local nature of the idea at the point $x$.  However, in several interesting cases in practice we can use this formulation to compute the mean deficit, because the tangent cone is isometric to a Euclidean space.  This occurs, for example, when the space $\X$ is a Hadamard manifold, in which case the tangent cone reduces to the traditional tangent space at $x$.  Another particularly interesting example for us is when the space $\X$ is a CAT(0) cubical complex and the point $x$ lies in the relative interior of a maximal cell. 

Mimicking precisely the proof of Proposition \ref{euclid-deficit}, we arrive at the following result.

\begin{thm}[Tangent criterion]
Consider a Hadamard space $\X$ and a non\-empty finite set \mbox{$A \subset \X$}.  At any point $x \not\in A$ at which the tangent cone $T_x \X$ is isometric to a Euclidean space, the mean deficit is the length of the shortest convex combination in $T_x \X$ of the vectors 
\[
\big(\big[[x,a]\big],d_a(x)\big) \qquad \mbox{for}~ a \in A.
\]
\end{thm}

\begin{cor}[Mean deficits in Hadamard manifolds]~ \label{deficit-manifold}
For any nonempty finite set $A$ in a Hadamard manifold $\X$, the mean deficit at any point $x \not\in A$  is the length of the shortest convex combination in the tangent space $T_x \X$ of the vectors 
$\mbox{\rm Log}_x(a)$ for $a \in A$.
\end{cor}

\section{Computing the mean set in Hadamard spaces} \label{sec:recognition}
Given a finite set $A$ in a Hadamard space $(\X,d)$, we seek to model its mean set computationally.  We introduce the following problems, for a given point $\bar x \in \X$: 
\newpage
\begin{quote}
\begin{description}
\item[decide]
whether or not $\bar x \in \mbox{mean}(A)$;
\item[certify]
 the property $\bar x \in \mbox{mean}(A)$ or the property $\bar x \not\in \mbox{mean}(A)$;
\item[compare]
$\bar x$ with other approximate mean points for $A$.
\end{description}
\end{quote}
Strictly speaking, we should refer to certifying $\bar x \not\in \mbox{minimal}(A)$, but in Hadamard space we know that the two sets $\mbox{mean}(A)$ and $\mbox{minimal}(A)$ coincide, so for simplicity we make no distinction.

The three questions are all straightforward in Hilbert space, as described in Proposition \ref{euclidean}.  Recognition and certification are also straightforward when the Hadamard space $\X$ is a manifold, as we outline in the following example.

\begin{exa}[Hadamard manifolds]
{\rm
Consider a finite subset $A$ of a Hadamard manifold $(\X,d)$, and a point $\bar x \in \X$.  In the tangent space $T_{\bar x}\X$, we first compute the differentials of the smooth maps $d_a^2$ at $\bar x$: 
\[
D(d_a^2)(\bar x) ~=~ -2\mbox{Log}_{\bar x}(a) ~\in~ T_{\bar x}\X.
\]
By convexity, the certificate property (\ref{certificate1}) is equivalent to
\bmye \label{certificate3}
w \in \Delta^A \qquad \mbox{and} \qquad \sum_{a \in A} w_a \mbox{Log}_{\bar x}(a) ~=~ 0.
\emye
Thus $\bar x$ lies in the mean set if and only if zero lies in the convex hull of the tangent vectors 
$\mbox{Log}_{\bar x}(a)$ (for $a \in A$), as we knew, less directly, from Corollary~\ref{deficit-manifold}:  the length of the shortest convex combination of the tangent vectors is the mean deficit.  One approach to the decision and certification questions is to check feasibility of  the system (\ref{certificate3}) via linear programming.  We either find a solution $w$ to the linear program, which satisfies the certificate property (\ref{certificate1}), or, via a separating hyperplane, find a tangent vector $z \in T_{\bar x}\X$ satisfying $\ip{\mbox{Log}_{\bar x}(a)}{z} > 0$ for all $a \in A$,
in which case the point $x = \mbox{Exp}_{\bar x}(tz)$ satisfies the certificate property~(\ref{certificate2}) for all small $t>0$.
}
\end{exa}

In more general Hadamard spaces, Proposition \ref{angle} allows us to recognize whether or not a point 
$\bar x$ lies in the mean set for a finite set $A$ via angles to geodesics between $\bar x$ and points in $A$.  We can therefore ``contract'' the set $A$ towards $\bar x$ without affecting the answer to the decision question, as described in the following result.

\begin{cor} \label{reduction}
Consider a finite subset $A$ of a Hadamard space $\X$ and a point $\bar x \not\in A$.  Corresponding to each point $a \in A$, consider any point $x_a$ on the geodesic segment $(\bar x,a]$.  Then $\bar x$ lies in the mean set for $A$ if and only if $\bar x$ lies in the mean set for the set $\{x_a : a \in A \}$.
\end{cor}

With the terminology of certificates, we can refine Corollary \ref{reduction} as follows.

\begin{thm}[Certifying means from contractions] \label{contraction}
Consider a finite subset $A$ of a Hadamard space $\X$ and a point $\bar x \not\in A$.  Corresponding to each point $a \in A$, suppose $0 < t_a \le 1$, and consider another subset of $\X$, defined by
\[
A' ~=~ \{t_a a + (1-t_a)\bar x : a \in A \}.
\]
If a vector $w \in \R^A$ certifies the property $\bar x \in \mbox{\rm mean}(A')$, then the vector with components $\frac{t_a w_a}{\sum_a t_a w_a }$ (for $a \in A$) certifies the property
$\bar x \in \mbox{\rm mean}(A)$.  On the other hand, if a point $x \in \X$ certifies the property
$\bar x \not\in \mbox{\rm mean}(A')$, then for all sufficiently small $t>0$, the point
$tx + (1-t)\bar x$ certifies the property $\bar x \not\in \mbox{\rm mean}(A)$.
\end{thm}

\pf
If the vector $w \in \Delta^A$ certifies the property $\bar x \in \mbox{\rm mean}(A')$, then by definition the point $\bar x$ is the corresponding barycenter for the set $A'$ , so Theorem \ref{sturm} and the CAT(0) inequality imply
\begin{eqnarray*}
0 
&\le& 
\sum_{a \in A} w_a 
\big( 
d^2(x,t_a a + (1-t_a)\bar x) - d^2(\bar x,t_a a + (1-t_a)\bar x) - d^2(x,\bar x)
\big) \\
&=& 
\sum_{a \in A} w_a 
\big( d^2(x,t_a a + (1-t_a)\bar x) - t_a^2 d^2(\bar x,a) - d^2(x,\bar x) \big) \\
&\le& 
\sum_{a \in A} w_a 
\big( 
t_a d^2(x,a) + (1\!-\!t_a)d(x,\bar x) - t_a(1\!-\!t_a)d^2(\bar x,a)
 - t_a^2 d^2(\bar x,a) - d^2(x,\bar x)
\big) \\
&=& 
\sum_{a \in A} w_a t_a
\big( 
d^2(x,a)  - d^2(\bar x,a) - d^2(x,\bar x)
\big).
\end{eqnarray*}
The first claim follows.

Now suppose that a point $x \in \X$ certifies the property 
$\bar x \not\in \mbox{\rm mean}(A')$.  For each point $a \in A$ define the corresponding point
$\hat a = t_a a + (1-t_a)\bar x$, so we have $d(x,\hat a) < d(\bar x,\hat a)$.
By the first variation formula, as $t \in (0,1]$ decreases to $0$ we have
\begin{eqnarray*}
0 ~>~ 
\frac{1}{t} \Big( d_{\hat a} \big( (1-t)\bar x + tx \big) - d_{\hat a}(\bar x) \Big) 
& \downarrow & 
-d(\bar x,x)\cos\angle x \bar x \hat a \qquad \mbox{and}  \\
\frac{1}{t} \Big( d_a \big( (1-t)\bar x + tx \big) - d_a(\bar x) \Big) 
& \downarrow & 
-d(\bar x,x)\cos\angle x \bar x a.
\end{eqnarray*}
The two right-hand sides are equal, so for all small $t>0$ we have our desired result:
$d_a \big( (1-t)\bar x + tx \big) < d_a(\bar x)$ for all $a \in A$.
\finpf

Our fundamental tool for deciding, with certification, whether or not a point $\bar x$ lies in the mean set for a finite set $A$ in a Hadamard space $\X$ is Proposition~\ref{recognition}:  $\bar x \in \mbox{mean}(A)$ if and only if it $\bar x$ minimizes the test function $f_{A,\bar x}$. Example \ref{pareto} shows that, if $\bar x \in \mbox{mean}(A)$, then that property is certified by an optimal solution of the dual problem~(\ref{p=2}), whereas if $\bar x \not\in \mbox{mean}(A)$, then that property is certified by any point $x$ satisfying $f_{A,\bar x}(x)<0$.  When $\X$ is a Hilbert space, Proposition \ref{euclidean} shows that the minimizer of the test function is even the nearest point in the mean set.

Typically, however, minimizing the test function does not produce the nearest point in the mean set, as the following example shows.

\begin{exa}[Nearest mean in a Hadamard space]
{\rm
In the Hadamard space $\{ x \in \R^2 : \min\{x_1,x_2\} \le 0\}$,
with the intrinsic metric induced by geodesics with respect to the Euclidean norm, consider the set
$A = \{ (0,1) \:,\: (\sqrt{3},0) \}$ and the point $\bar x = (0,-1)$.
The convex hull of $A$ is just the geodesic segment between its two points, which is the union of the segments $[(0,1) \:,\: (0,0) ]$  and $[ (0,0) \:,\: (\sqrt{3},0) ]$, and this coincides with its mean set.  The closest point to $\bar x$ on this geodesic is the point $(0,0)$.  The function $f$ is also minimized on this geodesic segment.  Its value at the point $(t,0)$, for $0 \le t \le \sqrt{3}$, is 
\[
\frac{1}{2} \max\{ (1+t)^2 - 4 \:,\: (\sqrt{3} - t)^2 - 4\},
\]
which is minimized by $t = \frac{1}{1+\sqrt{3}}$.  The corresponding point 
$(\frac{1}{1+\sqrt{3}},0)$ is not the closest mean to~$\bar x$.
}
\end{exa}

In general, the problem of estimating the distance in a Hadamard space $\X$ from a point $\bar x$ to the mean set for a finite subset $A$ does not seem easy. Certificates for the property $\bar x \not\in \mbox{mean}(A)$ do at least furnish easy lower bounds, as shown below.

\begin{prop}[Certified lower bound] \label{certified}
For a finite subset $A$ of a Hadamard space $(\X,d)$, if some point $x \in \X$ certifies the property
$\bar x \not\in \mbox{\rm mean}(A)$, then the distance from $\bar x$ to the mean set for $A$ is at least 
$\min_{a \in A} \{d_a(\bar x) - d_a(x) \}$.
\end{prop}

\pf
For any point $x' \in \X$ satisfying $d(\bar x,x') < d_a(\bar x) - d_a(x)$ for all points $a \in A$, we have
\[
d_a(x') ~\ge~ d_a(\bar x) - d(\bar x, x') ~>~ d_a(x)
\]
so $x'$ lies outside the mean set.
\finpf

\noindent
However, this lower bound on the distance to the mean set may never be tight, for any certificate.  Consider the example $A = \big\{ (0,-1) \:,\: (0,1) \big\}$, in $\R^2$, and $\bar x = (1,0)$.
The lower bound never exceeds $\sqrt{2} - 1$, but the distance to the mean set is $1$.

Trying to compute certificates or distance estimates by minimizing the test function is, in any case, not practical in general.  Variants of the proximal point method are valid in general Hadamard spaces \cite{bacak}, but computing proximal steps for the test function seems challenging.  More tractable is the dual problem
(\ref{p=2}), which involves maximizing the concave dual objective over the Euclidean simplex.  However, each dual objective evaluation requires a barycenter computation, for which general-purpose algorithms converge only very slowly. 

Since the problem of recognizing and certifying means in a general Hadamard space seems hard, we now narrow our focus to one specific class.  That class --- CAT(0) cubical complexes --- supports a more tractable approach to recognition and certification, and furthermore allows us to compare the relative accuracy of approximate mean points using the idea of mean discrepancy introduced in Section~\ref{mean-deficit}.

\section{Means in CAT(0) cubical complexes} \label{sec:interior}
A CAT(0) cubical complex is a particular kind of Hadamard space $(\X,d)$ that consists of a union of certain distinguished subsets called cells.  Each cell $C$ is the isometric image of a corresponding cube $\hat C$ with sides of unit length in some corresponding ambient Euclidean space:  informally, we identify $C$ and $\hat C$.  In particular, the isometry transfers the convex-analytic structure of $\hat C$ onto $C$, thereby defining the faces of $C$, and its relative interior, $\mbox{relint}(C)$.  Each face of a cell is also a cell.  Any two cells are either disjoint or are glued together along a unique common face.  A cell is {\em maximal} if it is not contained in any strictly larger cell.  The metric $d$ is just the intrinsic metric induced by the Euclidean distance on each cell.  The CAT(0) condition guaranteeing that a simply connected cubical complex is a Hadamard space is verifiable locally by checking a certain combinatorial condition at each vertex \cite{gromov}.  We start with a simple observation.

\begin{prop} \label{trio}
Consider any cell $C$ in a CAT(0) cubical complex $\X$, and any finite set $A \subset C$.  Denote by $\E$ an ambient Euclidean space for $C$.  Then the following three subsets of $C$ coincide.
\begin{itemize}
\item
The convex hull in the Euclidean cube $C \subset \E$ of the set $A \subset C$.
\item
The closed convex hull in the space $\X$ of the set $A \subset \X$.
\item
The mean set in the space $\X$ for the set $A \subset \X$.
\end{itemize}
Corresponding to any weight vector $w$ in the simplex $\Delta^A$, the barycenter of $A$ in 
$\X$ coincides with the barycenter of $A$ in $\E$, which is the convex combination $\sum_{a \in A} w_a a$.
\end{prop}

\pf
Denote the three sets, in order, by $\mbox{conv}_C(A)$, $\mbox{cl}\,\mbox{conv}(A)$, and $\mbox{mean}(A)$.  Corollary \ref{topology} implies 
$\mbox{mean}(A) \subset \mbox{cl}\,\mbox{conv}(A)$.  Since $\mbox{conv}_C(A)$ is a closed convex set in 
$\X$ that contains the set $A$, we deduce $\mbox{cl}\,\mbox{conv}(A) \subset \mbox{conv}_C(A)$.  To conclude that the three sets coincide, it remains to prove 
$\mbox{conv}_C(A) \subset \mbox{mean}(A)$.

Consider any weight vector $w$ in the simplex $\Delta^A$.  The corresponding barycenter of $A$ in the space $\E$ is the unique minimizer $\bar x$ of the function $f \colon \E \to \R$ defined by
\[
f(x) ~=~ \sum_{a \in A} w_a |x - a|^2 \qquad (x \in \E),
\]
namely the convex combination $\sum_{a \in A} w_a a$.
On the other hand, the corresponding barycenter of $A$ in the space $\X$ is the unique minimizer $\hat x$ of the function $\ip{w}{d_A^2}$, and satisfies $\hat x \in \mbox{mean}(A)$.  Since we know $\mbox{mean}(A) \subset C$, we know $\hat x \in C$, so
\[
\sum_{a \in A} w_a d_a^2(\hat x) ~=~ f(\hat x) ~\ge~ f(\bar x) ~=~ \sum_{a \in A} w_a d_a^2(\bar x).
\]
The uniqueness of $\hat x$ implies $\hat x = \bar x$, proving the final claim.

Any element of $\mbox{conv}_C(A)$ is a barycenter of $A$ in the space $\E$, by Proposition~\ref{euclidean}, and hence is also a barycenter of $A$ in the space $\X$, by the preceding argument, and hence lies in 
$\mbox{mean}(A)$.  We have proved $\mbox{conv}_C(A) \subset \mbox{mean}(A)$, as required.
\finpf

\begin{prop}[Relative interiors of maximal cells] \label{interior}
Consider any maximal cell $C$ in a CAT(0) cubical complex $\X$, any point $x$ in the relative interior of $C$.  Then for any point $a \in \X$, the open geodesic segment $(x,a)$ intersects $C$.
\end{prop}

\pf
As explained in \cite{ardila}, the geodesic $[x,a]$ decomposes into a chain of nontrivial geodesics 
$[x,y_1]$, $[y_1,y_2]$, \ldots. in corresponding cells $Q_1$, $Q_2$,\ldots. The cells $C$ and $Q_1$ must share a common face $F$ containing the point $x$.  Since $x$ lies in the relative interior of $C$, we deduce $F=C$.  Since the maximal cell $C$ is therefore a face of $Q_1$, we must have $C=Q_1$.  We deduce
$(x,y_1) \subset (x,a) \cap C$, completing the proof.
\finpf

Given a point $\bar x$ in the relative interior of a maximal cell $C$ in a CAT(0) cubical complex $\X$, we can now describe a simple procedure for recognizing whether or not $\bar x$ lies in the mean set for a finite set $A \subset \X$.  We first contract each point in $A$ towards $\bar x$, until each new point lies in the cube $C$.  We next reverse the contraction, but now in an ambient Euclidean space $\E$ for $C$, arriving at a new set $\hat A$.  By construction, the corresponding geodesics in $\X$ and in $\E$ have identical initial segments in $C$.  We can therefore solve the Euclidean recognition problem, by projecting $\bar x$ onto the convex hull of $\hat A$, and use the answer to that Euclidean problem to solve the original problem.  We describe the idea more formally in the following result.

\begin{thm}[Certifying relative interior means]\hspace{1pt} \label{certifying}
Consider any CAT(0) cubical complex $\X$, any finite set $A \subset \X$, and any point $\bar x \not\in A$. Suppose that $\bar x$ lies in the relative interior of a maximal cell $C$ with ambient Euclidean space $\E$.  For each point $a \in A$, if $t_a \in (0,1]$ is sufficiently small,
then the point $x_a = (1-t_a)\bar x + t_a a$ lies in $C$.
Given any such value $t_a$, define the point $z_a = \bar x + \frac{1}{t_a}(x_a - \bar x)$ in $\E$.  Then the mean deficit is given by
\bmye \label{straightened}
\mbox{\rm def}_A(\bar x) ~=~ \min \Big\{ |z - \bar x| : z =  \sum_{a \in A} w_a z_a, ~ z \in \E,~ w \in \Delta^A \Big\}.
\emye
Consider any pair $(\hat z,\hat w) \in \E \times \Delta^A$ attaining the minimum.  If the mean deficit is zero, then $\bar x \in \mbox{\rm mean}(A)$, and this property is certified by the vector $\hat w$.  On the other hand, if the mean deficit is strictly positive, then $\bar x \not\in \mbox{\rm mean}(A)$, and this property is certified by the point $(1-t)\bar x + t\hat z$ for all sufficiently small $t>0$.
\end{thm}

\pf
The fact that the point $x_a$ lies in the cube $C$ for all small $t_a > 0$ follows from Proposition \ref{interior}.  Fix a radius $\delta > 0$ such that the ball $\{y \in \X : d(\bar x,y) < \delta \}$ is contained in $C$.  For each point $y$ in this ball and each point $a \in A$, we have
\[
d_a(\bar x) = |\bar x-z_a| \qquad \mbox{and} \qquad \angle a\bar x y = \angle_\E \: z_a \bar x y.
\]
Consider the set in the Euclidean space $\E$ defined by $A' = \{ z_a : a \in A \}$.
By definition, we have
\begin{eqnarray*}
\mbox{\rm def}_A(\bar x)
&=&
\Big( \sup_{0 < d(\bar x,y) < \delta} \: \min_{a \in A} \{ d_a(\bar x)\cos\angle a\bar x y \} \Big)^+ \\
&=&
\Big( \sup_{0 < |\bar x-y|_\E < \delta} \: \min_{a \in A} \{ |\bar x-z_a|_\E \cos\angle_\E z_a \bar x y \} \Big)^+  ~=~
\mbox{\rm def}_{A'}(\bar x).
\end{eqnarray*}
Equation (\ref{straightened}) now follows by Proposition \ref{euclid-deficit}.

Suppose $\mbox{\rm def}_A(\bar x) = 0$.  We then know $\mbox{\rm def}_{A'}(\bar x) = 0$, so Proposition \ref{euclid-deficit} implies $\bar x \in \mbox{conv}(A')$, and hence
\[
\bar x ~=~ \sum_{a \in A} \hat w_a z_a ~=~ \sum_{a \in A} \hat w_a \Big(\bar x + \frac{1}{t_a}(x_a - \bar x)\Big).
\]
Consider the subset of $C$ defined by $A'' = \{ x_a : a \in A \}$.
The vector $w \in \Delta^A$ with components
\[
w_a ~=~ \frac{\frac{1}{t_a} \hat w_a}{\sum_{a'} \frac{1}{t_{a'}} \hat w_{a'} } \qquad (a \in A)
\]
certifies the property $\bar x \in \mbox{mean}(A'')$.  Theorem \ref{contraction} then implies that the vector $\hat w$ certifies the property $\bar x \in \mbox{mean}(A)$.

Suppose $\mbox{\rm def}_A(\bar x) > 0$.  We then know $\mbox{\rm def}_{A'}(\bar x) > 0$, so Proposition \ref{euclid-deficit} implies $\bar x \not\in \mbox{conv}(A')$.  By Proposition \ref{euclidean}, we deduce
\[
\ip{\hat z - \bar x}{x_a - \bar x} 
~ = ~
t_a \ip{\hat z - \bar x}{z_a - \bar x} ~>~ 0 \qquad \mbox{for all}~ a \in A,
\]
so for all small $t>0$ we have
\begin{eqnarray*}
0 & > & |t(\bar z - \bar x)|^2 + 2\ip{t(\bar z - \bar x)}{\bar x-x_a} ~=~
\big|\big((1-t)\bar x + t\hat z \big) - x_a \big|^2 ~-~ |\bar x-x_a|^2 \\
& = & d_{x_a}^2\big((1-t)\bar x + t\hat z \big) - d_{x_a}^2(\bar x).
\end{eqnarray*}
This shows that the point $(1-t)\bar x + t\hat z$ certifies the property $\bar x \not\in \mbox{mean}(A')$, so the final claim follows by Theorem \ref{contraction}.
\finpf

Theorem \ref{certifying} validates the following procedure.  In the complex $\X$, the procedure only involves computing distances and midpoints;  in the ambient Euclidean space for the cube $C$, the procedure computes affine combinations of pairs of points, along with a single projection problem onto a polytope at each iteration.

\begin{alg}[Decide if $\bar x \in \mbox{mean}(A)$ in CAT(0) cubical complex $\X$] \label{algorithm}
{\rm
\begin{algorithmic}
\STATE
\STATE 	{\bf input:} nonempty finite $A \subset \X$, maximal cell $C$ in $\X$, \\
ambient Euclidean space $\E$ for $C$ \\
point $\bar x \in \mbox{relint}(C) \setminus A$, tolerance $\epsilon \ge 0$
\FOR{$a \in A$}
\STATE	$d_a = d(a,\bar x)$
\STATE	$x_a = a$									\hfill \% {\em initial point in direction of $a$}
\STATE	$t_a = 1$									\hfill \% {\em initial proportion of distance to $a$}
\WHILE{$x_a \not\in C$}
\STATE	$x_a = \frac{1}{2}x_a + \frac{1}{2}\bar x$	\hfill \% {\em in $\X$, move $x_a$ towards $\bar x$}
\STATE	$t_a = \frac{1}{2} t_a$			\hfill \% {\em new proportion:  $x_a = (1-t_a)\bar x + t_a a$}
\ENDWHILE
\STATE	in $\E$, define $z_a = \bar x + \frac{1}{t_a}(x_a - \bar x)$
\hfill \% {\em so  $x_a = (1-t_a)\bar x + t_a z_a$ in $\E$}
\ENDFOR
\STATE	choose $(z,w)$ attaining 			\hfill \% {\em in $\E$, project $\bar x$ onto 
$\mbox{\rm conv}\{z_a\}$}
\[
\mbox{mean$\underline~$deficit} ~=~ \min_{z \in \E,\, w \in \Delta^A} \Big\{ |z - \bar x| : z = \sum_{a \in A} w_a z_a \Big\}
\]
\IF{$\mbox{mean$\underline~$deficit} > \epsilon$}
\FOR{$a \in A$}
\WHILE{$z \not\in C$ \OR $d(z,a) \ge d_a$}
\STATE	in $\E$ define $z = \frac{1}{2} z + \frac{1}{2}\bar x$ \hfill \% {\em contract into $C$ and closer to all $a \in A$}
\ENDWHILE
\ENDFOR
\RETURN{``The point $z$ certifies that the point $\bar x$ lies outside the mean set.''}
\ELSE
\RETURN{``The point $\bar x$ is an approximate barycenter, with weight vector $w$.''}
\ENDIF
\end{algorithmic}
}
\end{alg}
\bigskip

By convexity, for all points $a \in A$, the quantity
\[
\frac{1}{t} \Big( d\big( (1-t)\bar x + ta , a \big) - d(\bar x,a) \Big)
\]
is a nondecreasing function of $t \in (0,1]$.  In the case 
``$\mbox{mean$\underline~$deficit} > \epsilon$'', after completing the {\tt for} loop, we have $d(z,a) < d(\bar x,a)$ for all $a \in A$, as required.

The mean deficit can help us explore the mean set computationally.  We begin with a simple two-dimensional illustration.

\begin{exa}[A mean set in a planar complex] \label{planar}
{\rm
Consider the CAT(0) cubical complex consisting of three squares in $\R^2$:
\[
[0,1] \times [-1,0] \hspace{15mm} [-1,0] \times [-1,0] \hspace{15mm} [-1,0] \times [0,1].
\]
Simply connected square complexes like this trivially satisfy the combinatorial condition \cite{gromov} equivalent to the CAT(0) property.  In this space, consider the set $A$ consisting of the three points
\[
a = [1 ~~ 0] \hspace{15mm} b = [-1 ~~ 0] \hspace{15mm} c = [0 ~~ 1].
\]
The mean set and convex hull of $A$ coincide, both consisting of the convex hull of the Euclidean triangle $\triangle 0bc$ along with the line segment $[0,a]$.  We can verify this observation computationally by repeated random trials, as follows.  For each trial, we randomly choose one of the three squares, $C$, sample a random point $\bar x$ from a uniform distribution on $C$, and then plot $\bar x$ with a shade determined by its mean deficit, calculated as in Algorithm \ref{algorithm}:  lighter shades indicate smaller mean deficits.  The first panel of Figure \ref{fig:squares1-squares2} is a typical illustration:  the mean set is clearly apparent, approximated by the lightly shaded region.

\begin{figure}
\centering
\parbox{6cm}
{
\begin{center}
\includegraphics[width=0.25\textwidth]{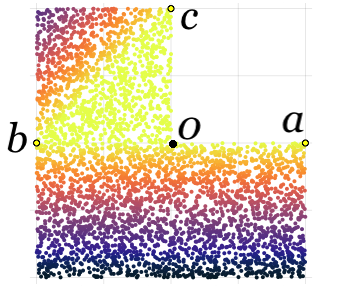}
\end{center}
}
\hspace{20mm}
\parbox{6cm}
{
\begin{center}
\includegraphics[width=0.4\textwidth]{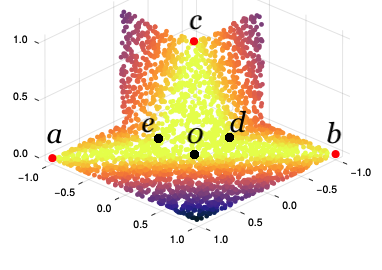}
\end{center}
}
\caption{Heat maps for the mean deficit of three-point sets $\{a,b,c\}$ in CAT(0) square complexes, illustrating Examples \ref{planar} and \ref{5squares} respectively.
}
\label{fig:squares1-squares2}
\end{figure}
}
\end{exa}

The next example involves a CAT(0) square complex embedded in $\R^3$.

\begin{exa}[Three points in a square complex] \label{5squares}
{\rm
Consider the CAT(0) cubical complex consisting of five squares in $\R^3$:
\begin{eqnarray*}
[0,1] \times [-1,0] \times \{0\} \hspace{15mm} \mbox{} & [0,1] \times [0,1] \times \{0\} 
& \mbox{} \hspace{15mm} [-1,0] \times [0,1] \times \{0\} \\
\mbox{$[-1,0]$} \times \{0\} \times [0,1]  & & \{0\} \times [-1,0] \times [0,1].
\end{eqnarray*}
In this space, consider the set $A$ consisting of the three points
\[
a = [1 ~~ -\!1 ~~ 0] \hspace{15mm} b = [-1 ~~ 1 ~~ 0] \hspace{15mm} c = [0 ~~ 0 ~~ 1].
\]
The geodesic $[a,b]$ passes through the point $0 = [0 ~~ 0 ~~ 0]$;  the geodesic $[b,c]$ passes through the point $d = [-\frac{1}{2} ~~ 0 ~~ 0]$; the geodesic $[c.a]$ passes through the point 
$e = [0 ~~ -\frac{1}{2} ~~ 0]$.  One can then verify formally that the mean set and convex hull of $A$ coincide, both consisting of the convex hulls of the four Euclidean triangles
\[
\triangle 0ae \qquad \triangle 0ec \qquad \triangle 0cd \qquad \triangle 0db.
\]
Much the same process as described in Example \ref{planar} yields the mean deficit heat map in the second panel of Figure \ref{fig:squares1-squares2}.  Once again, the mean set is clearly apparent, approximated by the lightly shaded region.
}
\end{exa}

\section{A mean set in a CAT(0) cubical complex} \label{mean-example}
\begin{figure}
\centering
\parbox{6cm}
{
\begin{center}
\includegraphics[width=0.4\textwidth]{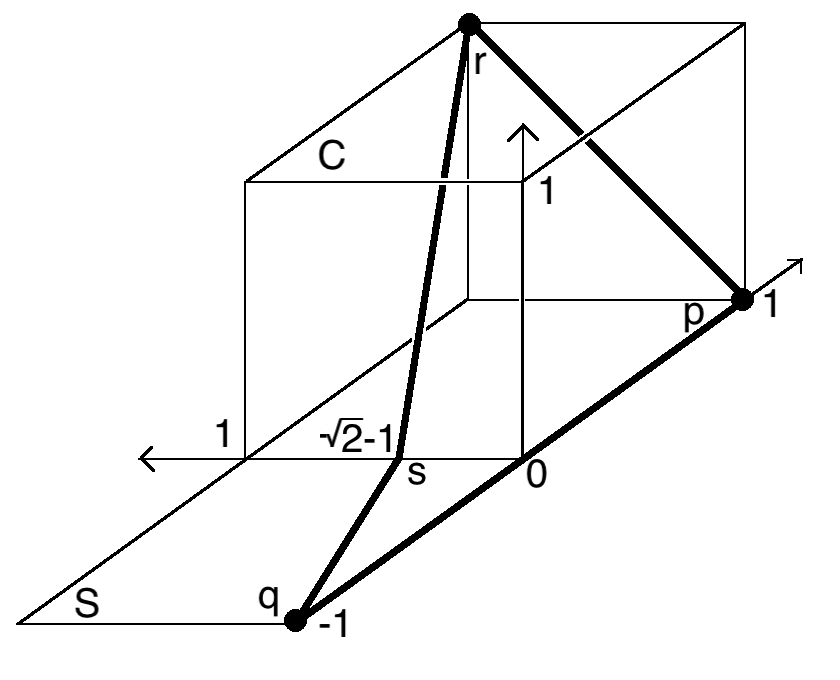}
\end{center}
}
\hspace{15mm}
\parbox{6cm}
{
\begin{center}
\includegraphics[width=0.4\textwidth]{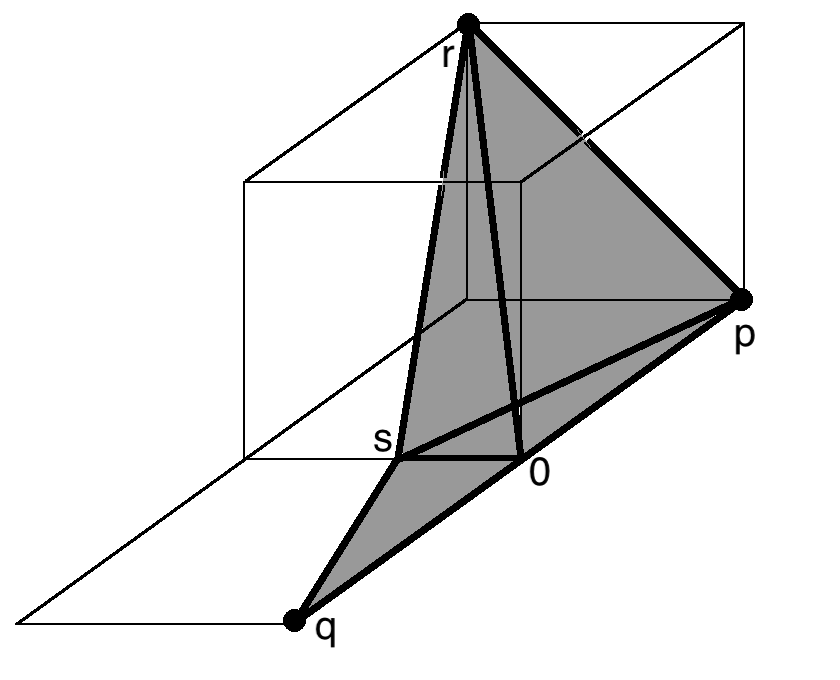}
\end{center}
}
\caption{The first panel shows a geodesic triangle with vertices $p,q,r$ in a CAT(0) cubical complex.  The second panel shows the convex hull of $p,q,r$, which contains the three-dimensional simplex with vertices $p,0,s,r$.}
\label{fig:complex3-complex4}
\end{figure}

We next consider in detail a simple example to illustrate the distinction between the mean set and the convex hull.  We consider a complex $\X$ consisting of two maximal cells --- a square $S$ and a three-dimensional cube $C$ --- glued together along a common edge $E$.  It is easy to verify that this complex is CAT(0).  We can embed $\X$ isometrically in the ambient Euclidean space $\R^3$.  For example, we can make the identification 
\[
C = [0,1]^3, \qquad S ~=~ [-1,0] \times [0,1] \times \{0\}, \qquad E ~=~ \{0\} \times [0,1] \times \{0\}. 
\]
In $\X$, we consider a set $A$ consisting of three points:
\[
p = [1~0~0], \qquad q = [-1~0~0], \qquad r = [1~1~1].
\]
The three geodesics comprising the corresponding triangle include two Euclidean line segments: $[p,q]$ (which passes through the point $0 = [0~0~0]$ and $[r,p]$).  The third geodesic, $[q,r]$, is the union of two Euclidean line segments:  $[q,s]$ in $S$ and $[s,r]$ in $C$, for a point $s$ on the edge $E$.  It is easy to verify $s = [0 ~ (\sqrt{2} - 1) ~ 0]$.
We illustrate the complex and triangle in the first panel of Figure \ref{fig:complex3-complex4}.

The convex hull of the set $A$ contains the points $0$ and $s$, and hence must also contain the Euclidean triangle $0qs$ in the square $S$ and the Euclidean simplex with vertices $0$, $p$, $r$, and $s$ in the cube $C$.  It is not hard to verify that the union of this triangle and simplex is convex in the complex $\X$, and hence is the (closed) convex hull of $A$.  Worth noting is that the convex hull of the three points in $A$, illustrated in the second panel of Figure \ref{fig:complex3-complex4}, contains a three-dimensional Euclidean ball.

The mean set for $A$ is a subset of the convex hull, by Corollary \ref{topology}, but is strictly smaller.  To identify it, we can appeal to Theorem \ref{certifying}.  That result shows in particular that, in the square $S$, any point in the relative interior of the Euclidean triangle $0qs$ must belong to $\mbox{mean}(A)$.  Since $\mbox{mean}(A)$ is closed, it must therefore contain the whole triangle $0qs$.

\begin{figure}
\centering
\parbox{6cm}
{
\begin{center}
\includegraphics[width=0.4\textwidth]{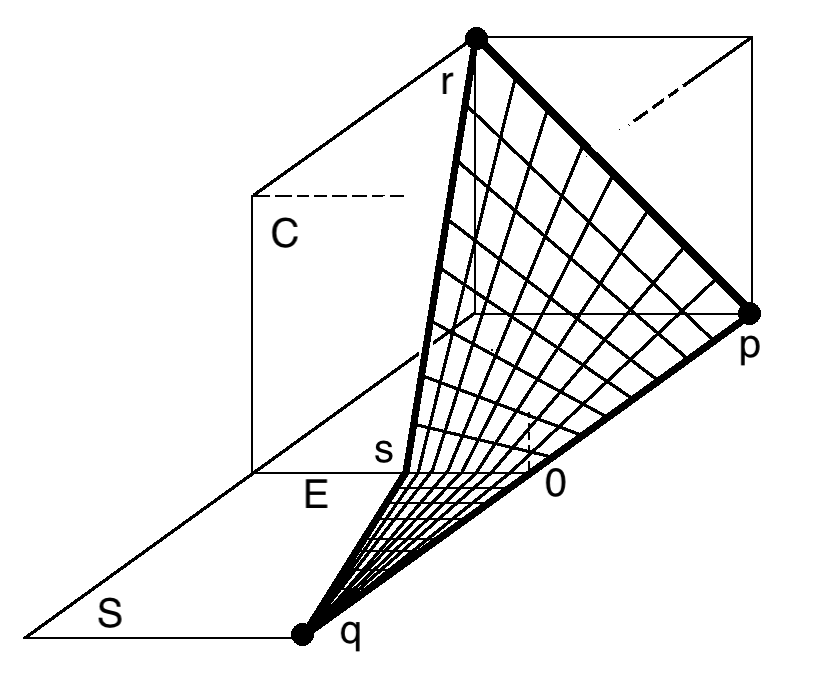}
\end{center}
}
\hspace{15mm}
\parbox{6cm}
{
\begin{center}
\includegraphics[width=0.4\textwidth]{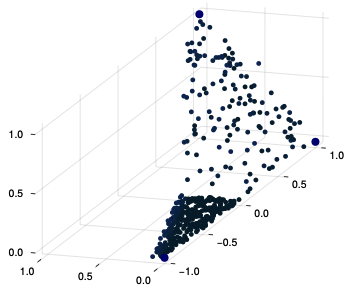}
\end{center}
}
\caption{The first panel shows the mean set of three points, $p,q,r$, in a CAT(0) cubical complex in $\R^3$, a two-dimensional semi-algebraic surface.  The second panel plots those randomly sampled points in the complex with mean deficit less than $\frac{1}{10}$.}
\label{fig:complex5-cube}
\end{figure}

Consider next any point $c = [x~y~z]$ in the relative interior of the cube $C$.  The geodesic $[c,q]$ crosses the common edge $E$ at the point 
\[
g ~=~ \Big[ 0 \quad \frac{y}{1+\sqrt{x^2 + z^2}} \quad 0 \Big].
\]
By Theorem \ref{certifying}, $c \in \mbox{mean}(A)$ holds if and only if $c$ is a convex combination of the points $p,r,g$ in the cube $C$.  In that case, the vectors $p-c$, $r-c$, and $g-c$ must be coplanar, so
\begin{eqnarray*}
0 &=& \det
\left[
\begin{array}{ccc}
x-1 & y & z \\
x-1 & y-1 & z-1 \\
x & y -  \frac{y}{1+\sqrt{x^2 + z^2}} & z
\end{array}
\right]
~=~
\det
\left[
\begin{array}{ccc}
x-1 & y & z \\
0 & -1 & -1 \\
1 & -  \frac{y}{1+\sqrt{x^2 + z^2}} & 0
\end{array}
\right] \\ \\
&=&
z-y + \frac{y(1-x)}{1+\sqrt{x^2 + z^2}}.
\end{eqnarray*}
For $x>0$, solving gives $y$ as a function of $x$ and $z$:
\bmye \label{surface}
y ~=~ \frac{z(1+\sqrt{x^2 + z^2})}{x+\sqrt{x^2 + z^2}}.
\emye
Each point $[x~y~z]$ in the relative interior of the simplex with vertices $p,0,s,r$ satisfies
$0<z<x<1$.  On the domain $\{(x,z) : 0<z<x<1\}$, equation (\ref{surface}) describes an algebraic function, whose graph, a two-dimensional algebraic surface, we can now verify is contained in 
$\mbox{mean}(A)$.  Its closure is the intersection of $\mbox{mean}(A)$ with $C$.

We sketch the mean set for $A$ in the first panel of Figure \ref{fig:complex5-cube}.
To be concrete, we consider some specific points.  The points
\[
\Big[4t \quad \frac{1}{3}(1+5t) \quad 3t \Big] \qquad \mbox{for}~ 0 < t < \frac{1}{4}
\]
all lie in the graph of the function (\ref{surface}) and hence in $\mbox{mean}(A)$.  On the other hand, the point $[\frac{1}{2} ~ \frac{1}{3} ~ \frac{1}{4}]$, which lies in the convex hull of $A$, does not satisfy equation (\ref{surface}) so lies outside $\mbox{mean}(A)$.

The second panel of Figure \ref{fig:complex5-cube} explores the same mean set using Algorithm \ref{algorithm}.  Specifically, we sample random points $\bar x$ from uniform distributions on the square $S$ and the cube $C$, and run Algorithm \ref{algorithm} with tolerance $\epsilon = \frac{1}{10}$ to decide whether or not to plot $\bar x$ is an approximate weighted mean.

\section{Cell boundaries in CAT(0) cubical complexes} \label{boundaries}
We again consider a finite set $A$ in a CAT(0) cubical complex $\X$.  We seek to decide whether or not a point $\bar x \in \X$ belongs to the mean set for $A$, and to certify the answer.  Section \ref{sec:interior} described a procedure for resolving this question when $\bar x$ lies in the relative interior of a maximal cell.  Here we answer the general question, using semidefinite programming. To describe the algorithm, we first observe how to construct semidefinite representations for directional derivatives of distance functions.

In the complex $\X$, consider any cell $C$ containing the point $\bar x$.  We can suppose that $C$ is embedded in a Euclidean space $\E_C$.  Given any point $a \in \X$, define a convex function 
$d_{a,C} \colon \E_C \to \R$ by
\[
d_{a,C}(y) ~=~ 
\left\{
\begin{array}{ll}
d_a(y) & (y \in C) \\
+\infty & (y \not\in C).
\end{array}
\right.
\]
We can compute the geodesic $[\bar x,a]$ in polynomial time~\cite{hayashi}, expressing it as a finite sequence of nontrivial segments, each confined in a single cell in the complex.  Denote the initial segment by $[\bar x, x_a]$, and the cell containing it by $Q_a$.  

For any direction $u$ in the tangent cone $T_C(\bar x)$ (a polyhedral cone in $\E_C$), the directional derivative $d'_{a,C}(\bar x;u)$ is the rate of increase of the distance function $d_a$ as we move from the point $\bar x$ in the cube $C$ in the direction $u$.  By the first Variation formula, that rate of change is the product of the norm of $u$ with the negative cosine of the angle at $\bar x$ that the geodesic $[\bar x,\bar x + \epsilon u]$ (for any small $\epsilon > 0$) makes with the geodesic $[\bar x,a]$, or equivalently, with the geodesic $[\bar x,x_a]$.  Consequently we have
\[
d'_{a,C}(\bar x;u) ~=~ d'_{x_a,C}(\bar x;u) \qquad \mbox{for all}~ u \in \E_C,
\]
both sides being $+\infty$ for $u \not\in T_C(\bar x)$.

The cells $C$ and $Q_a$ have a unique common face $F_a$, and for any point $y$ in the cell $C$, since the geodesic $[x_a,y]$ must pass through $F_a$, we have
\[
d(x_a,y) ~=~ \min_{x \in F_a} \{ |x_a -x|_{Q_a} + |x -y|_C \},
\]
where $|\cdot|_{Q_a}$ and $|\cdot|_C$ denote the Euclidean distances in the cells $Q_a$ and $C$ respectively.  This formula, via standard tools \cite{mod_opt}, gives a semidefinite representation of the function $d_{x_a,C}$, and hence also of its subdifferential 
$\partial d_{x_a,C}(\bar x) \subset \E_C$,
from which we deduce a semidefinite representation of its support function, namely
\[
d'_{x_a,C}(\bar x;\cdot) \colon \E_C \to \bar\R.
\]
We conclude that we can efficiently construct a semidefinite representation of the directional derivative function $d'_{a,C}(\bar x;\cdot) \colon \E_C \to \bar\R$.  Denote this function by $f_{a,C}$.

We now return to the problem of recognizing whether the point $\bar x$ lies in the mean set for a given finite set of points $A \subset \X$.  In the complex, it suffices to restrict attention to maximal cells $C$ containing $\bar x$.  For each cube $C$ in the set ${\mathcal C}$ of maximal such cells, we first solve the optimization problem
\[
\mbox{($\mbox{P}_C$)} \qquad
\inf_{u \in \E_C} \max_{a \in A} f_{a,C}(u).
\]
Our observation above allows us to frame ($\mbox{P}_C$) as an explicit semidefinite program, which we can then solve efficiently.  By positive homogeneity, its optimal value is either $0$ or $-\infty$.  In the latter case, we find a vector $u \in \E_C$ satisfying $f_{a,C}(u) < 0$ for all points $a \in A$.  By bisection search, we can find a number $\epsilon > 0$ such that $d_a(\bar x + \epsilon u) < d_a(\bar x)$ for all $a \in A$, certifying that $\bar x$ is outside the mean set for $A$.

Suppose, on the other have, that we find each problem ($\mbox{P}_C$), for cells $C \in {\mathcal C}$, has optimal value zero.  Consider our usual test function $f \colon \X \to \R$, defined by 
\[
f(x) ~=~ \frac{1}{2} \max_{a \in A} \{ d^2_a(x) - d^2_a(\bar x) \}.
\]
By assumption, the point $\bar x$ minimizes $f$ over each cell $C$, since otherwise, some point $x_0 \in C$ satisfies $d_a(x_0) < d_a(\bar x)$ for all points $a \in A$, and then by convexity we deduce the contradiction $f_{a,C}(x_0 - \bar x) < 0$ for all $a \in A$.
Thus $\bar x$ locally minimizes the test function, which is convex, so in fact
$\bar x$ is a global minimizer, proving that $\bar x$ lies in the mean set for $A$.  It remains only to find a certificate for that fact.

We seek a certificate $w$ in the simplex $\Delta^A$ such that the point $\bar x$ minimizes the function $\ip{w}{d^2_A}$ over the space $\X$.  We know that such a certificate exists, by Theorem \ref{coincide}.  We can suppose $\bar x \not\in A$, since otherwise the certificate is obvious.

By convexity, $w \in \Delta^A$ is a certificate if and only if, for all cells $C \in {\mathcal C}$, we have
\[
0 ~\le~ \Big( \sum_{a \in A} w_a d^2_{a,C} \Big)'(\bar x;u) ~=~ 
2\sum_{a \in A} w_a d_a(\bar x) d'_{a,C}(\bar x;u) \qquad \mbox{for all}~ u \in \E_C,
\]
which, by positive homogeneity, we can rewrite as
\[
\inf_{u \in \E_C} \sum_{a \in A} w_a d_a(\bar x) f_{a,C}(u) ~=~ 0 \qquad \mbox{for all}~ C \in {\mathcal C}.
\]
Equivalently, we seek a vector $v$ in the intersection of the sets 
\bmye \label{inf}
V_C ~=~ 
\Big\{ v \in \Delta^A : \inf_{u \in \E_C} \sum_{a \in A} v_a f_{a,C}(u) ~=~ 0 \Big\} 
\qquad (C \in {\mathcal C}).
\emye
Given such a vector, we set $w_a = \frac{v_a}{d_a(\bar x)}$ and rescale to get a certificate $w \in \Delta^A$.

For each cube $C \in {\mathcal C}$, define a function $g_C \colon \R^A \to \bar\R$ by
\[
g_C(z) ~=~ \inf_{u \in \E_C} \max_{a \in A} \{ f_{a,C}(u) + z_a \}.
\]
Since each function $f_{a,C}$ is semidefinite representable, so is the function $g_C$.  Notice $g_C(0) = 0$.
We claim $V_C = \partial g_C(0)$.

To prove this claim, notice $v \in \partial g_C(0)$ if and only if
\bmye \label{subdifferential}
\ip{v}{z} ~\le~ \max_{a \in A}  \{ f_{a,C}(u) + z_a \} \qquad \mbox{for all}~ z \in \R^A,~ u \in \E_C.
\emye
Suppose that property (\ref{subdifferential}) holds.
Setting $u=0$ implies $w \in \Delta^A$.  Setting $z_a = -f_{a,C}(u)$ for all points $a \in A$, and then $u=0$, implies $v \in V_C$.  Conversely, $v \in V_C$ implies, for all $u \in \E_C$ and $z \in \R^A$, the inequalities
\[
\max_{a \in A}  \{ f_{a,C}(u) + z_a \} ~\ge~ \sum_{a \in A} v_a(f_{a,C}(u) + z_a) ~\ge~ \ip{v}{z},
\]
so property (\ref{subdifferential}) holds, as we claimed.

The functions $g_C$, for $C \in {\mathcal C}$, are all semidefinite representable, and hence so are their subdifferentials $\partial g_C(0) = V_C$.  We can therefore express the problem of finding a vector in the (nonempty) intersection 
$\cap_C V_C$, as a semidefinite program, which we can solve efficiently, and immediately deduce a certificate.

\begin{figure}
\centering
\begin{center}
\includegraphics[width=60mm]{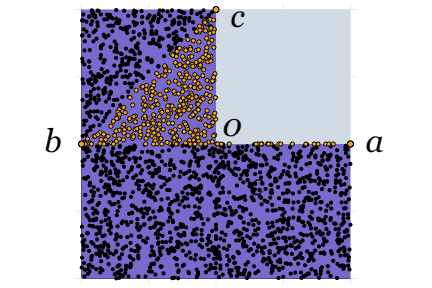}
\end{center}
\caption{A scatter plot of the mean set for three points $a,b,c$ in a CAT(0) square complex, testing cell boundary points by semidefinite programming.}
\label{fig:3squares}
\end{figure}

To illustrate the computational technique described in this section, we return to Example \ref{planar}.  As before, we explore the mean set by repeated random trials:  in each trial, we randomly choose one of the three squares, $C$, sample a random point $\bar x$ from a uniform distribution on C, solve the semidefinite program ($\mbox{P}_C$) to determine whether or not $\bar x$ lies in the mean set, and, in Figure \ref{fig:3squares}, plot it lightly shaded if so, and darkly shaded otherwise.  However, unlike the heat map in Figure \ref{fig:squares1-squares2}, the geodesic segment $[0,a]$ on the relative boundary of a square would not be apparent in the mean set using these samples alone:  $[0,a]$ almost surely contains none of the sample points.  To generate the scatter plot in Figure \ref{fig:3squares},  we include some random trial points on the line segment $[a,b]$, all of which appear in the final plot.

\bibliographystyle{plain}
\small
\parsep 0pt

\begin{thebibliography}{10}

\bibitem{ardila}
F.~Ardila, M.~Owen, and S.~Sullivant.
\newblock Geodesics in {CAT(0)} cubical complexes.
\newblock {\em Advances in Applied Mathematics}, 48:142--163, 2012.

\bibitem{arnaudon-nielsen}
M.~Arnaudon and F.~Nielsen.
\newblock On approximating the {R}iemannian $1$-center.
\newblock {\em Computational Geometry}, 46:93--104, 2013.

\bibitem{basso}
G.~Basso, Y.~Krifka, and E.~Soultanis.
\newblock A non-compact convex hull in generalized non-positive curvature.
\newblock {\tt arXiv:2301.03835}, 2023.

\bibitem{bacak}
M.~Ba\u{c}\'{a}k.
\newblock The proximal point algorithm in metric spaces.
\newblock {\em Isr. J. Math.}, 194:689--701, 2013.

\bibitem{bacak-convex}
M.~Ba\u{c}\'{a}k.
\newblock {\em Convex Analysis and Optimization in Hadamard Spaces}.
\newblock De Gruyter, Berlin, 2014.

\bibitem{bacak-open}
M.~Ba\u{c}\'{a}k.
\newblock Old and new challenges in {H}adamard spaces.
\newblock {\em Japanese Journal of Mathematics}, 18:115--168, 2023.

\bibitem{mod_opt}
A.~Ben-Tal and A.~Nemirovski.
\newblock {\em Lectures on Modern Convex Optimization}.
\newblock SIAM, Philadelphia, 2001.

\bibitem{bessenyi}
M.~Bessenyi, D.~Cs. Kert\'esz, and R.L. Lovas.
\newblock A sandwich with segment convexity.
\newblock {\em Journal of Mathematical Analysis and Applications}, 500, 2021.

\bibitem{billera}
L.J. Billera, S.P. Holmes, and K.~Vogtmann.
\newblock Geometry of the space of phylogenetic trees.
\newblock {\em Advances in Applied Mathematics}, 27:733--767, 2001.

\bibitem{bridson}
M.R. Bridson and A.~Haefliger.
\newblock {\em Metric Spaces of Non-Positive Curvature}.
\newblock Springer-Verlag, Berlin, 1999.

\bibitem{burago}
D.~Burago, Yu~Burago, and S.~Ivanov.
\newblock {\em A Course in Metric Geometry}.
\newblock Amer. Math. Soc., Providence, RI, 2001.

\bibitem{cartan}
E.~Cartan.
\newblock {\em Le\c{c}ons sur la {G}\'eom\'etrie des {E}spaces de {R}iemann}.
\newblock Gauthier-Villars, Paris, 1928.

\bibitem{descombes}
D.~Descombes and U.~Lang.
\newblock Convex geodesic bicombings and hyperbolicity.
\newblock {\em Geometriae Dedicata}, 177:367--384, 2015.

\bibitem{espinola}
R.~Esp\'inola and A.~Nicolae.
\newblock Geodesic {P}tolemy spaces and fixed points.
\newblock {\em Nonlinear Analysis}, 74:27--34, 2011.

\bibitem{gromov}
M.~Gromov.
\newblock {\em Structures m\'etriques pour les vari\'et\'es {R}iemanniennes},
  volume~1 of {\em Textes Math\'ematiques.}
\newblock CEDIC, Paris, Edited by J. Lafontaine and P. Pansu, 1981.

\bibitem{gromov99}
M.~Gromov.
\newblock Metric structures for {R}iemannian and non-{Ri}emannian spaces.
\newblock In {\em Progress in Mathematics}, volume 152 (Based on a 1981 French
  original.). Birkh\"{a}user Boston Inc, Boston MA, 1999.

\bibitem{hansen}
P.~Hansen, M.~Labb\'e, D.~Peeters, and J.-F. Thisse.
\newblock Single facility location on networks.
\newblock {\em Annals of Discrete Mathematics}, 31:113--146, 1987.

\bibitem{hayashi}
K.~Hayashi.
\newblock A polynomial time algorithm to compute geodesics in {CAT(0)} cubical
  complexes.
\newblock {\em Discrete and Computational Geometry}, 65:636--654, 2021.

\bibitem{jost94}
J.~Jost.
\newblock Equilibrium maps between metric spaces.
\newblock {\em Calculus of Variations and Partial Differential Equations},
  2:173--204, 1994.

\bibitem{lewis-lopez-nicolae-basic}
A.S. Lewis, G.~L\'opez, and A.~Nicolae.
\newblock Basic convex analysis in metric spaces with bounded curvature.
\newblock {\em SIAM Journal on Optimization}, 34:366--388, 2023.

\bibitem{lubiw}
A.~Lubiw, D.~Maftuleac, and M.~Owen.
\newblock Shortest paths and convex hulls in 2{D} complexes with non-positive
  curvature.
\newblock {\em Computational Geometry}, 89, 2020.

\bibitem{lytchak-generic}
A.~Lytchak and A.~Petrunin.
\newblock About every convex set in any generic {R}iemannian manifold.
\newblock {\em Journal f\"{u}r die reine und angewandte Mathematik},
  782:235--245, 2022.

\bibitem{papadopoulos}
A.~Papadopoulos.
\newblock {\em Metric Spaces, Convexity and Nonpositive Curvature}.
\newblock European Mathematical Society, Z{\"u}rich, 2005.

\bibitem{sturm2}
K.-T. Sturm.
\newblock Probability measures on metric spaces of nonpositive curvature.
\newblock {\em Contemporary Mathematics}, 338:357--390, 2003.

\bibitem{takahashi}
W.~Takahashi.
\newblock A convexity in metric space and nonexpansive mappings, {I}.
\newblock {\em Kodai Mathematical Seminar Reports}, 22:142--149, 1970.

\end{thebibliography}

\def\cprime{$'$} \def\cprime{$'$}

\end{document}